\documentclass[aop,MSNbibl,dvips]{arximspdf}
\usepackage{mathbh}
\usepackage{graphicx}

\doi{10.1214/12-AOP749} %kopijuoti is PTS
\volume{41}
\issue{3A}
\pubyear{2013}
\firstpage{1427}
\lastpage{1460}

\makeatletter

\newcommand{\rrvert}{\vert}
\newcommand{\llvert}{\vert}
\newcommand{\mathds}{\mathbh}
\newcommand{\iint}{\int\!\!\!\int}
\newcommand{\eqref}[1]{(\ref{#1})}

\newtheorem{thmm}{Theorem}[section]
\newtheorem{cor}[thmm]{Corollary}
\newtheorem{lem}[thmm]{Lemma}
\newtheorem{prop}[thmm]{Proposition}
\newproclaim{defn}[thmm]{Definition}
\newproclaim{ex}{Example}[section]
\newproclaim{rem}[thmm]{Remark}
\newproclaim{Conclusion}{Conclusion}
\newtheorem{thmmmm}{Theorem}[section]
\newtheorem{propp}[thmmmm]{Proposition}
\newtheorem{lemm}[thmmmm]{Lemma}
\newproclaim{remm}[thmmmm]{Remark}
\makeatother

\begin{document}
\begin{frontmatter}

\title{Convergence to the equilibria for self-stabilizing processes in
double-well landscape\thanksref{T1}}
\runtitle{Convergence for self-stabilizing processes}
\thankstext{T1}{Supported by the DFG-funded CRC 701, Spectral
Structures and Topological Methods in
Mathematics, at the University of Bielefeld.}

\begin{aug}
\author{\fnms{Julian} \snm{Tugaut}\corref{}\ead[label=e1]{jtugaut@math.uni-bielefeld.de}}
\runauthor{J. Tugaut}
\affiliation{Universit\"at Bielefeld}
\address{Fakult\"at f\"ur Mathematik\\
Universit\"at Bielefeld\\
D-33615 Bielefeld\\
Germany\\
\printead{e1}}
\end{aug}

% HISTORY:
\received{\smonth{12} \syear{2010}}
\revised{\smonth{12} \syear{2011}}

% ABSTRACT
%
\begin{abstract}
We investigate the convergence of McKean--Vlasov diffusions in a
nonconvex landscape. These processes are linked to nonlinear partial
differential equations. According to our previous results, there are at
least three stationary measures under simple assumptions. Hence, the
convergence problem is not classical like in the convex case. By using
the method in Benedetto et al. [\textit{J. Statist. Phys.} \textbf{91} (1998)
1261--1271]
about the monotonicity of the free-energy, and
combining this with a complete description of the set of the stationary
measures, we prove the global convergence of the self-stabilizing processes.
\end{abstract}

% KEYWORDS
%


\end{frontmatter}
\section*{Introduction}\label{intr}
We investigate the weak convergence in long-time of the following
so-called self-stabilizing process:
%e1 #&#
%
\renewcommand{\theequation}{\Roman{equation}}
\begin{equation}
\label{eqintroinit} %\tag{I}
\cases{ %
\displaystyle X_t=X_0+\sqrt{\varepsilon}B_t-\int
_0^tV' (X_s )\,ds-\int
_0^tF'\ast u_s^\varepsilon
(X_s )\,ds,
\vspace*{2pt}\cr
u_s^\varepsilon=\mathcal{L} (X_s ) .}
\end{equation}
Here, $\ast$ denotes the convolution. Since the own law of the process
intervenes in the drift, this equation is nonlinear, in the sense of
McKean. We note that $X_t$ depends on $\varepsilon$. We do not write
$\varepsilon$ for simplifying the reading.

The motion of the process is generated by three concurrent forces. The
first one is the derivative of a potential $V$---the
confining potential. The second influence is a Brownian motion
$ (B_t )_{t\in\mathbb{R}_+}$.
It allows the particle to move upwards the potential $V$. The third
term---the so-called self-stabilizing term---represents
the attraction between all the others trajectories. Indeed, we remark:
$F'\ast u_s^\varepsilon (X_s(\omega_0) )=\int_{\omega\in\Omega
}F' (X_s(\omega_0)-X_s(\omega) )\,d\mathbb{P} (\omega )$
where $ (\Omega,\mathcal{F},\mathbb{P} )$ is the underlying
measurable space.

This kind of processes were introduced by McKean, see~\cite{McKean} or
\cite{McKean1966}. Here, we will make some
smoothness assumptions on the interaction potential~$F$. Let just note
that it is possible to consider nonsmooth $F$. If $F$ is the Heaviside
step function and $V:=0$, \eqref{eqintroinit} is the Burgers
equation; see~\cite{SV1979}. If $F:=\delta_0$, and without confining
potential, it is the Oelschl\"ager equation, studied in~\cite{Oel1985}.

The particle $X_t$ which verifies \eqref{eqintroinit} can be seen as
one particle in a continuous mean-field system of an infinite number of
particles. The mean-field system that we will consider is a random
dynamical system like
%e2 #&#
%
\begin{equation}
\label{eqintromeanfield} %\tag{II}
\cases{ %
\displaystyle dX_t^1=\sqrt{\varepsilon}\,dB_t^1-V'
\bigl(X_t^1 \bigr)\,dt-\frac{1}{N}\sum
_{j=1}^NF' \bigl(X_t^1-X_t^j
\bigr)\,dt,
\cr
\vdots
\cr
\displaystyle dX_t^i=\sqrt{\varepsilon}\,dB_t^i-V'
\bigl(X_t^i \bigr)\,dt-\frac{1}{N}\sum
_{j=1}^NF' \bigl(X_t^i-X_t^j
\bigr)\,dt,
\cr
\vdots
\cr
\displaystyle dX_t^N=\sqrt{\varepsilon}\,dB_t^N-V'
\bigl(X_t^N \bigr)\,dt-\frac{1}{N}\sum
_{j=1}^NF' \bigl(X_t^N-X_t^j
\bigr)\,dt, }
\end{equation}
where the $N$ Brownian motions $ (B_t^i )_{t\in\mathbb{R}_+}$
are independents. Mean-field systems are the subject of a rich
literature; see~\cite{DG1987} for the large deviations for $N\to+\infty
$ and~\cite{M1996} under weak assumptions on $V$ and $F$. For
applications, see~\cite{CDPS2010} for social interactions or \cite
{CX2010} for the stochastic partial differential equations.

The link between the self-stabilizing process and the mean-field system
when $N$ goes to $+\infty$ is called the propagation of chaos; see \cite
{Sznitman} under Lipschitz properties;~\cite{BRTV} if $V$ is a constant;
\cite{Malrieu2001} or~\cite{Malrieu2003} when both potentials are
convex;~\cite{BAZ1999} for a more precise result;~\cite{BGV2007,DPPH1996} or~\cite{DG1987} for a sharp estimate;~\cite{CGM} for a
uniform result in time in the nonuniformly convex case.

Equation \eqref{eqintromeanfield} can be rewritten in the following way:
%e3 #&#
%
\setcounter{equation}{1}
\begin{equation}
 d\mathcal{X}_t=\sqrt{\varepsilon}
\mathcal{B}_t-N\nabla\Upsilon^N (\mathcal{X}_t
)\,dt,
\end{equation}
where the $i$th coordinate of $\mathcal{X}_t$ (resp., $\mathcal{B}_t$)
is $X_t^i$ (resp., $B_t^i$) and
\[
\Upsilon^N (\mathcal{X} ):=\frac{1}{N}\sum
_{j=1}^NV (\mathcal{X}_j )+
\frac{1}{2N^2}\sum_{i=1}^N\sum
_{j=1}^NF (\mathcal{X}_i-
\mathcal{X}_j )
\]
for all $\mathcal{X}\in\mathbb{R}^N$. As noted in~\cite{TT}, the
potential $\Upsilon^N$ converges toward a functional~$\Upsilon$ acting
on the measures. A perturbation (proportional to $\varepsilon$) of
$\Upsilon$ will play the central role in the article.

As observed in~\cite{DG1987}, the empirical law of the mean-field
system can be seen as a perturbation of the law of the diffusion \eqref
{eqintroinit}. Consequently, the long-time behavior of $\mathcal
{L} (X_t )$ that we study in this paper provides some
consequences on the exit time for the particle system \eqref
{eqintromeanfield}.

Also, the convergence plays an important role in the exit problem for
the self-stabilizing process since the exit time is strongly linked to
the drift according to the Kramers law (see~\cite{DZ} or~\cite{HIP})
which converges toward a homogeneous function if the law of the
process converges toward a stationary measure.

Let us recall briefly some of the previous results on diffusions like
\eqref{eqintroinit}. The existence problem has been investigated by
two different methods. The first one consists in the application of a
fixed point theorem; see~\cite{McKean,BRTV,CGM} or \cite
{HIP} in the nonconvex case. The other consists of a propagation of
chaos; see, for example,~\cite{M1996}. Moreover, it has been proved in
Theorem 2.13 in~\cite{HIP} that there is a unique strong solution.

In~\cite{McKean}, the author proved---by using Weyl's lemma---that the
law of the unique strong solution $du_t^\varepsilon$ admits a $\mathcal
{C}^\infty$-continuous density $u_t^\varepsilon$ with respect to the
Lebesgue measure for all $t>0$. Furthermore, this density satisfies a
nonlinear partial differential equation of the following type:
%e4 #&#
%
\begin{equation}
\label{eqintropde} %\tag{III}
\frac{\partial}{\partial t}u_t^\varepsilon(x)=
\frac{\partial}{\partial
x} \biggl\{\frac{\varepsilon}{2}\frac{\partial}{\partial x}u_t^\varepsilon
(x)+u_t^\varepsilon(x) \bigl(V'(x)+F'
\ast u_t^\varepsilon(x) \bigr) \biggr\} .
\end{equation}
It is then possible to study equations like \eqref{eqintropde} by
probabilistic methods which involve diffusions \eqref{eqintroinit} or
\eqref{eqintromeanfield}; see~\cite{CGM,Funaki1984,Malrieu2003}.
Reciprocally, equation \eqref{eqintropde} is a useful tool for
characterizing the stationary measure(s) and the long-time behavior;
see~\cite{BRTV,BRV,Tamura194,Tamura1987,Veret2006}. In~\cite{HT1}, in
the nonconvex case, by using \eqref{eqintropde}, it has been proved
that the diffusion \eqref{eqintroinit} admits at least three
stationary measures under assumptions easy to verify. One is symmetric,
and the two others are not. Moreover, Theorem~3.2 in~\cite{HT1} states
the thirdness of the stationary measures if $V''$ is convex and $F'$ is
linear. This nonuniqueness prevents the long-time behavior from being
as intuitive as in the case of unique stationary measure.

The work in~\cite{HT2} and~\cite{HT3} provides some estimates of the
small-noise asymptotic of these three stationary measures. In
particular, the convergence toward Dirac measures and its rate of
convergence have been investigated. This will be one of the two main
tools for obtaining the convergence.

Convergence for \eqref{eqintroinit} is not a new subject. In \cite
{BRV}, if $V$ is identically equal to~$0$, the authors proved the
convergence toward the stationary measure by using an
ultracontractivity property, a Poincar\'e inequality and a comparison
lemma for \mbox{stochastic} processes. The ultracontractivity property still
holds if $V$ is not convex by using the results in~\cite{KKR}. It is
possible to conserve the Poincar\'e inequality by using the theorem of
Muckenhoupt (see~\cite{logsob2000}) instead of the Bakry--Emery
theorem. But, the comparison lemma needs some convexity properties.
However, it is possible to apply these results\vadjust{\goodbreak} if the initial law is
symmetric in the synchronized case ($V''(0)+F''(0)\geq0$); see Theorem~7.10 in~\cite{TT}.

Another method consists of using the propagation of chaos in order to
derive the convergence of the self-stabilizing process from the one of
the mean-field system. However, we shall use it independently of the
time and the classical result which is on a finite interval of time is
not sufficiently strong. Cattiaux, Guillin and Malrieu proceeded a
uniform propagation of chaos in~\cite{CGM} and obtained the convergence
in the convex case, including the nonuniformly strictly convex case.
See also~\cite{Malrieu2003}. Nevertheless, according to Proposition
5.17 and Remark 5.18 in~\cite{TT}, it is impossible to find a general
result of uniform propagation of chaos. In the synchronized case, if
the initial law is symmetric, it is possible to find such a uniform
propagation of chaos; see Theorems 7.11 and 7.12 in~\cite{TT}.

The method that we will use in this paper is based on the one of \cite
{BCCP}. See also \cite
{Malrieu2003,Tamura194,Malrieu2001,HS1987,AMTU2001} for the convex
case. In the nonconvex case, Carrillo, McCann and Villani provide the
convergence in~\cite{CMV2003} under two restrictions: the center of
mass is fixed and $V''(0)+F''(0)>0$ (that means it is the synchronized
case).

However, by combining the results in~\cite{HT1,HT2,HT3} with the work
of~\cite{BCCP} (and the more rigorous proofs in~\cite{CMV2003} about
the free-energy), we will be able to prove the convergence in a more
general setting. The principal tool of the paper is the monotonicity of
the free-energy along the trajectories of \eqref{eqintropde}.

First, we introduce the following functional:
%e5 #&#
%
\begin{equation}
\label{eqintrometapotentiel} %\tag{IV}
\Upsilon(u):=\int_\mathbb{R}
V(x)\,du(x)+\frac{1}{2}\iint_{\mathbb
{R}^2}F(x-y)\,du(x)\,du(y) .
\end{equation}
This quantity appears intuitively as the limit of the potential in
\eqref{eqintromeanfield} for $N\to+\infty$. We consider now the
free-energy of the self-stabilizing process~\eqref{eqintroinit},
\[
\Upsilon_\varepsilon(u):=\frac{\varepsilon}{2}\int_\mathbb{R}
u(x)\log \bigl(u(x)\bigr)\,dx+\Upsilon(u)
\]
for all measures $du$ which are absolutely continuous with respect to
the Lebesgue measure. We can note that $du_t^\varepsilon$ satisfies this
hypothesis for all $t>0$.

The paper is organized as follows. After presenting the assumptions, we
will state the first results, in particular, the convergence of a
subsequence $ (u_{t_k}^\varepsilon )_k$. This subconvergence will
be used for improving the results about the thirdness of the stationary
measures. Then, we will give the main statement which is the
convergence toward a stationary measure, briefly discuss the
assumptions of the theorem and give the proof. Subsequently, we will
study the basins of attraction by two different methods and prove that
these basins are not reduced to a single point. Finally, we postpone
four results in the annex, including Proposition~\ref{propafm} which
extends the classical higher-bound for the moments of the
self-stabilizing processes.

\subsection*{Assumptions}
\begin{figure}

\includegraphics{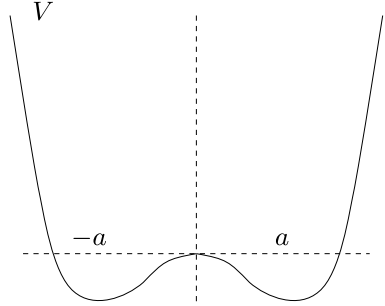}

\caption{Potential $V$.}\label{fig1}
\end{figure}

We assume the following properties on the confining potential~$V$ (see Figure~\ref{fig1}):
\begin{longlist}[(V-1)]
\item[(V-1)] $V$ is an even polynomial function with $\deg(V)=:2m\geq4$.
\item[(V-2)] The equation $V'(x)=0$ admits exactly three solutions:
$a$, $-a$ and $0$ with $a>0$. Furthermore, $V''(a)>0$ and $V''(0)<0$.
Then, the bottoms of the wells are located in $x=a$ and $x=-a$.
\item[(V-3)] $V(x)\geq C_4x^4-C_2x^2$ for all $x\in\mathbb{R}$ with $C_2,C_4>0$.
\item[(V-4)] $ \lim_{x\to\pm\infty}V''(x)=+\infty$ and
$V''(x)>0$ for all $x\geq a$.
\item[(V-5)] $V''$ is convex.
\item[(V-6)] Initialization: $V(0)=0$.
\end{longlist}
%
%

%
%%\begin{figure}[H]
%% \psfrag{$V$}{$V$}
%% \psfrag{$-a$}{$-a$}
%% \psfrag{$a$}{$a$}
%% \centerline{\includegraphics[width=4.5cm]{poten.eps}}
%%
%% \end{figure}
%%

%f1 #&#
%

Let us remark that the positivity of $V''$ on $[-a;a]^c$ \mbox{[in hypothesis
(V-4)]} is an immediate consequence of (V-1) and (V-5). The simplest and
most studied example is $V(x):=\frac{x^4}{4}-\frac{x^2}{2}$. Also, we
would like to stress that weaker assumptions could be considered, but
all the mathematical difficulties are present in the polynomial case,
and it allows us to avoid some technical and tedious computations. Let
us present now the assumptions on the interaction potential $F$:
\begin{longlist}[(F-1)]
\item[(F-1)] $F$ is an even polynomial function with $\deg(F)=:2n\geq2$.
\item[(F-2)] $F$ and $F''$ are convex.
\item[(F-3)] Initialization: $F(0)=0$.
\end{longlist}
Under these assumptions, we know by~\cite{HT1} that \eqref
{eqintroinit} admits at least one symmetric stationary measure. And,
if $\sum_{p=0}^{2n-2}\frac{\llvert F^{(p+2)}(a)\rrvert }{p!}a^p<F''(0)+V''(a)$, there are at least two asymmetric stationary
measures: $u^\varepsilon_+$ and $u^\varepsilon_-$. Furthermore, we know by
\cite{HT2} that there is a unique nonnegative real $x_0$ such that
$V'(x_0)+\frac{1}{2}F'(2x_0)=0$ and $V''(x_0)+\frac
{F''(0)+F''(2x_0)}{2}>0$. The same paper provides that $u^\varepsilon_0$
converges weakly toward $\frac{1}{2}\delta_{x_0}+\frac{1}{2}\delta_{-x_0}$ and $u^\varepsilon_\pm$ converges weakly toward $\delta_{\pm a}$
in the small-noise limit.

We present now the assumptions on the initial law $du_0$:
\begin{longlist}[(ES)]
\item[(ES)] The $8q^2$th moment of the measure $du_0$ is finite with
$q:=\max \{m,n \}$.
\item[(FE)] The probability measure $du_0$ admits a $\mathcal{C}^\infty
$-continuous density $u_0$ with respect to the Lebesgue measure. And,
the entropy $\int_\mathbb{R} u_0\log(u_0)$ is finite.
\end{longlist}
Under (ES), \eqref{eqintroinit} admits a unique strong solution.
Indeed, the assumptions of Theorem 2.13 in~\cite{HIP} are satisfied:
$V'$ and $F'$ are locally Lipschitz, $F'$ is odd, $F'$ grows
polynomially, $V'$ is continuously differentiable, and there exists a
compact $K$ such that $V''$ is uniformly negative on $K^c$. Moreover,
we have the following inequality:
%e6 #&#
%
\begin{equation}
\label{majoration} %\tag{V}
\max_{1\leq j\leq8q^2}\sup_{t\in\mathbb{R}_+}\mathbb{E}
\bigl[\llvert X_t\rrvert^{j} \bigr]\leq M_0 .
\end{equation}
We deduce immediately that the family $ (u_t^\varepsilon )_{t\in
\mathbb{R}_+}$ is tight. The assumption (FE) ensures that the initial
free-energy is finite. In the following, we shall use occasionally one
of the following three additional properties concerning the two
potentials $V$ and $F$ and the initial law $du_0$:
\begin{longlist}[(SYN)]
\item[(LIN)] $F'$ is linear.
\item[(SYN)] $V''(0)+F''(0)>0$.
\item[(FM)] For all $N\in\mathbb{N}$, we have $\int_\mathbb{R}
|x|^N\,du_0(x)<+\infty$.
\end{longlist}
In the following, three important properties linked to the enumeration
of the stationary measures for the self-stabilizing process \eqref
{eqintroinit} will be helpful for proving the convergence:
\begin{longlist}[(0M1)]
\item[(M3)] The process \eqref{eqintroinit} admits exactly three
stationary measures. One is symmetric: $u^\varepsilon_0$ and the other
ones are asymmetric: $u^\varepsilon_+$ and $u^\varepsilon_-$. Furthermore,
$\Upsilon_\varepsilon(u^\varepsilon_+)=\Upsilon_\varepsilon(u^\varepsilon_-)<\Upsilon_\varepsilon(u^\varepsilon_0)$.
\item[(M3)$'$] There exists $M>0$ such that the diffusion \eqref
{eqintroinit} admits exactly three stationary measures with
free-energy less than $M$. Furthermore, we have $\Upsilon_\varepsilon
(u^\varepsilon_+)=\Upsilon_\varepsilon(u^\varepsilon_-)<\Upsilon_\varepsilon
(u^\varepsilon_0)\leq M$ ; $u^\varepsilon_0$ is symmetric, and $u^\varepsilon_+$
and $u^\varepsilon_-$ are asymmetric.
\item[(0M1)] The process \eqref{eqintroinit} admits only one
symmetric stationary measure $u^\varepsilon_0$.
\end{longlist}
In the following, we will give some simple conditions such that (M3),
(M3)$'$ or (0M1) are true.

Finally, we recall assumption (H) introduced in~\cite{HT2}:
\begin{longlist}
\item[(H)] A family of measures $ (v^\varepsilon )_\varepsilon$
verifies assumption (H) if the family of positive reals $ (\int_\mathbb{R} x^{2n}v^\varepsilon(dx) )_{\varepsilon>0}$ is bounded.
\end{longlist}
The aim of the weaker assumption (M3)$'$ is to obtain the convergence
even if there exists a family of stationary measures which does not
verify the assumption (H).

For concluding the \hyperref[intr]{Introduction}, we write the statement of the main
theorem:

\begin{thmmm*}
Let $du_0$ be a probability measure which verifies
\textup{(FE)} and \textup{(FM)}. Under \textup{(M3)}, $u_t^\varepsilon$ converges weakly toward a
stationary measure.\vspace*{-2pt}
\end{thmmm*}

%s1 #&#
\section{First results}
This section is devoted to present the tools that we will use for
proving the main result of the paper. Furthermore, we provide some new
results about the thirdness of the stationary measures for the
self-stabilizing processes.

We introduce the following functional:
\[
\Upsilon_\varepsilon^-(u):=\frac{\varepsilon}{2}\int_\mathbb{R}
u(x)\log \bigl(u(x) \bigr)\mathds{1}_{\{u(x)<1\}}\,dx+\int_\mathbb{R}
V(x)u(x)\,dx .
\]
This new functional is linked to the free-energy $\Upsilon_\varepsilon$.
The interaction part and the positive contribution of the entropy term
$\int_\mathbb{R} u\log (u )$ have been removed. Let us
consider a measure $u$ which verifies the previous assumptions. Due to
the nonnegativity of the functions $F$ and $x\mapsto\frac{\varepsilon
}{2}u(x)\log (u(x) )\mathds{1}_{\{u(x)\geq1\}}$, we obtain
directly the inequality $\Upsilon_\varepsilon(u)\geq\Upsilon_\varepsilon^-(u)$.

In the following, we will need two particular functions [the
free-energy of the system and a function $\eta_t$ such that $\frac
{d}{dt}u_t^\varepsilon(x)=\frac{d}{dx}\eta_t(x)$].\vspace*{-2pt}
%de1.1 #&#
%
\begin{defn}
\label{defintroxi}
For all $t\in\mathbb{R}_+$, we introduce the functions
\[
\xi(t):=\Upsilon_\varepsilon \bigl(u_t^\varepsilon \bigr)
\quad\mbox{and}\quad  \eta_t:=\frac{\varepsilon}{2}\frac{\partial}{\partial x}u_t^\varepsilon
+u_t^\varepsilon \bigl(V'+F'\ast
u_t^\varepsilon \bigr) .\vspace*{-2pt}
\]
\end{defn}
According to \eqref{eqintropde}, we remark that if $\eta_t$ is
identically equal to $0$, then $u_t^\varepsilon$ is a~stationary measure
for \eqref{eqintroinit}.

We recall the following well-known entropy dissipation:\vspace*{-2pt}
%pr1.2 #&#
%
\begin{prop}
\label{propintrodecrease}
Let $du_0$ be a probability measure which verifies \textup{(FE)} and~\textup{(ES)}. Then,
for all $t,s\geq0$, we have
\[
\xi(t+s)\leq\xi(t)\leq\xi(0)<+\infty .
\]
Furthermore, $\xi$ is derivable, and we have
\[
\xi'(t)\leq-\int_\mathbb{R}\frac{\eta_t^2}{u_t^\varepsilon} .\vspace*{-2pt}
\]
\end{prop}
See~\cite{CMV2003} for a proof.\vspace*{-2pt}

%s1.1 #&#
\subsection{Preliminaries}
Let us introduce the functional space
\[
\mathcal{M}_{8q^2}:= \biggl\{f\in\mathcal{C}_0^2
(\mathbb{R},\mathbb {R}_+ )   \Big|  \int_\mathbb{R} f(x)\,dx=1
\biggr\} .
\]
We can remark that $u_t^\varepsilon\in\mathcal{M}_{8q^2}$ for all $t>0$;
see~\cite{McKean}. The first tool is the Proposition \ref
{propintrodecrease} [i.e., to say the fact that the free-energy is
decreasing along the orbits of~\eqref{eqintropde}]. The second one is
its lower-bound.\vspace*{-2pt}
%le1.3 #&#
%
\begin{lem}
\label{lemfrminoration}
$\!\!\!$There exists $\Xi_varepsilon\in\mathbb{R}$ such that \mbox{$
\inf_{u\in\mathcal{M}_{8q^2}}\Upsilon_\varepsilon(u)\geq\Xi_\varepsilon$}.\vadjust{\goodbreak}
\end{lem}
\begin{pf}
Let us recall $\Upsilon_\varepsilon(u)\geq\Upsilon_\varepsilon^-(u)$. It
suffices then to prove the inequality $\inf_{u\in\mathcal
{M}_{8q^2}}\Upsilon_\varepsilon^-(u)\geq\Xi_\varepsilon$. We proceed as in
the first part of the proof of Theorem 2.1 in~\cite{BCCP}. We show that
we can minorate the negative part of the entropy by a function of the
second moment. Then a growth condition of $V$ will provide the result.

We split the negative part of the entropy into two integrals,
\[
-\int_\mathbb{R} u(x)\log \bigl(u(x) \bigr)
\mathds{1}_{\{u(x)<1\}
}\,dx=-I_+-I_-\vspace*{-2pt}
\]
with
\[
I_+:=\int_\mathbb{R} u(x)\log \bigl(u(x)
\bigr)\mathds {1}_{\{e^{-|x|}<u(x)<1\}}\,dx\vspace*{-2pt}
\]
and
\[ I_-:=\int_\mathbb{R} u(x)\log \bigl(u(x)
\bigr)\mathds {1}_{\{u(x)\leq e^{-|x|}\}}\,dx .\vspace*{-2pt}
\]
By definition of $I_+$, we have the following estimate:
\begin{eqnarray*}
I_+&\geq&\int_\mathbb{R} u(x)\log \bigl(e^{-|x|} \bigr)
\mathds{1}_{\{
e^{-|x|}<u(x)<1\}}\,dx
\\[-3pt]
&\geq&-\int_\mathbb{R} |x|u(x)\mathds{1}_{\{e^{-|x|}<u(x)<1\}}\,dx
\\[-3pt]
&\geq&-\int_\mathbb{R} |x|u(x)\,dx\geq-\frac{1}{2}-
\frac{1}{2}\int_\mathbb {R} x^2u(x)\,dx .
\end{eqnarray*}
By putting $\gamma(x):=\sqrt{x}\log(x)\mathds{1}_{\{x<1\}}$, a simple
computation provides $\gamma(x)\geq-2e^{-1}$ for all $x<1$. We deduce
\begin{eqnarray*}
I_-=\int_\mathbb{R}\sqrt{u(x)}\gamma\bigl(u(x)\bigr)
\mathds{1}_{\{u(x)\leq
e^{-|x|}\}}\,dx\geq-2e^{-1}\int_\mathbb{R}
e^{-{|x|}/{2}}\,dx=-8e^{-1} .
\end{eqnarray*}
Consequently, it yields
\[
-\int_\mathbb{R} u(x)\log \bigl(u(x) \bigr)
\mathds{1}_{\{u(x)<1\}}\,dx\leq \frac{1}{2}\int_\mathbb{R}
x^2u(x)\,dx+\frac{1}{2}+8e^{-1} .
\]
This implies
%e7 #&#
%
\renewcommand{\theequation}{\arabic{section}.\arabic{equation}}
\setcounter{equation}{0}
\begin{equation}
\label{eqfrminoration} \Upsilon_\varepsilon^-(u)\geq-\frac{\varepsilon}{4}-4
\varepsilon e^{-1}+\int_\mathbb{R} \biggl(V(x)-
\frac{\varepsilon}{4}x^2 \biggr)u(x)\,dx .
\end{equation}
By hypothesis, there exist $C_2, C_4>0$ such that $V(x)\geq
C_4x^4-C_2x^2$ so the function $x\mapsto V(x)-\frac{\varepsilon}{4}x^2$ is
lower-bounded by a negative constant. This achieves the proof.\vspace*{-2pt}
\end{pf}
Let us note that the unique assumption we used is $ \lim_{x\to\pm\infty}V''(x)=+\infty$.\vspace*{-2pt}
%le1.4 #&#

\begin{lem}
\label{lemfrconv}
Let $du_0$ be a probability measure which satisfies the assumptions
\textup{(FE)} and \textup{(ES)}. Then, there exists $L_0\in\mathbb{R}$ such that
$\Upsilon_\varepsilon (u_t^\varepsilon )$ converges toward~$L_0$ as time goes
to infinity.\vadjust{\goodbreak}
\end{lem}
\begin{pf}
The assumption (FE) implies $\xi(0)=\Upsilon_\varepsilon(u_0)<\infty$. As
$\xi$ is nonincreasing by Lemma~\ref{propintrodecrease} and
lower-bounded by a constant $\Xi_\varepsilon$ according to Lemma~\ref
{lemfrminoration}, we deduce that the function $\xi$ converges
toward a real $L_0$.\vspace*{-2pt}
\end{pf}
%
%le1.5 #&#
%
\begin{lem}
\label{lemfrstation}
If and only if $\xi'(t)=0$, the following is true: $u_t^\varepsilon$ is a
stationary measure $u^\varepsilon$.\vspace*{-2pt}
\end{lem}
\begin{pf}
If $u_t^\varepsilon$ is a stationary measure $u^\varepsilon$, then $\xi
(t)=\Upsilon_\varepsilon (u_t^\varepsilon )=\Upsilon_\varepsilon
(u^\varepsilon)$ is a constant. This provides $\xi'(t)=0$.

Reciprocally, if $\xi'(t)=0$, Proposition~\ref{propintrodecrease} implies
\[
\int_\mathbb{R}\frac{\eta_t^2}{u_t^\varepsilon}=0 .
\]
We deduce $\eta_t(x)=0$ for all $x\in\mathbb{R}$. This means that
$u_t^\varepsilon$ is a stationary measure.\vspace*{-2pt}
\end{pf}
%
%s1.2 #&#
\subsection{Subconvergence}
%th1 #&#
%
\begin{thmm}
\label{thmfrsubconv}
Let $du_0$ be a probability measure which satisfies the assumptions
\textup{(FE)} and \textup{(ES)}. Then there exists a stationary measure $u^\varepsilon$ and
a sequence $(t_k)_{k}$ which converges toward infinity such that
$u_{t_k}^\varepsilon$ converges weakly toward $u^\varepsilon$.\vspace*{-2pt}
\end{thmm}
\begin{pf}
\textit{Plan}: First, we use the convergence of $\int_t^\infty\xi'(s)\,ds$
toward $0$ when $t$ goes to infinity, and we deduce the existence of a
sequence $(t_k)_k$ such that $\xi' (t_k )$ tends toward $0$
when $k$ goes to infinity. Then, we extract a subsequence of $(t_k)_k$
for obtaining an adherence value. By using a test function, we prove
that this adherence value is a stationary measure.
\begin{longlist}[\textit{Step} 1.]
\item[\textit{Step} 1.] Lemma~\ref{lemfrconv} implies that $\int_t^\infty\xi'(s)\,ds$
collapses at infinity. According to Proposition \ref
{propintrodecrease}, the sign of $\xi'$ is a constant, so we deduce
the existence of an increasing sequence $ (t_k )_{k\in\mathbb
{N}}$ which goes to infinity such that $\xi'(t_k)\longrightarrow0$.

\item[\textit{Step} 2.] The uniform boundedness of the first $8q^2$ moments with
respect to the time allows us to use Prohorov's theorem: we can extract
a subsequence [we continue to write it $ (t_k )_{k}$ for
simplifying] such that $u_{t_k}^\varepsilon$ converges weakly toward a
probability measure $u^\varepsilon$.

\item[\textit{Step} 3.] We consider now a function $\varphi\in\mathcal{C}^\infty
(\mathbb{R},\mathbb{R} )\cap\mathcal{L}_2 (u^\varepsilon
)$ with compact support, and we estimate the following quantity:
\begin{eqnarray*}
&&\biggl\llvert \int_\mathbb{R}\varphi(x) \biggl\{
\frac{\varepsilon}{2}\frac{\partial
}{\partial x}u_{t_k}^\varepsilon(x)+u_{t_k}^\varepsilon(x)
\bigl[V'(x)+\bigl(F'\ast u_{t_k}^\varepsilon
\bigr) (x) \bigr] \biggr\}\,dx\biggr\rrvert
\\[-3pt]
&&\qquad=\biggl\llvert \int_\mathbb{R}\varphi(x)\eta_{t_k}(x)\,dx
\biggr\rrvert =\biggl\llvert \int_\mathbb{R}\varphi(x)\sqrt
{u_{t_k}^\varepsilon(x)}\frac{\llvert \eta_{t_k}(x)\rrvert }{\sqrt{u_{t_k}^\varepsilon(x)}}\,dx\biggr\rrvert
\\
&&\qquad\leq \biggl(\int_\mathbb{R}\varphi(x)^2u_{t_k}^\varepsilon(x)\,dx
\biggr)^{{1}/{2}}\times \biggl(\int_\mathbb{R}
\frac{1}{u_{t_k}^\varepsilon
(x)} \bigl(\eta_{t_k}(x) \bigr)^2\,dx
\biggr)^{{1}/{2}}
\\
&&\qquad\leq\sqrt{-\xi'(t_k)}\sqrt{\int
_\mathbb{R}\varphi(x)^2u_{t_k}^\varepsilon
(x)}\longrightarrow0
\end{eqnarray*}
when $k$ goes to infinity; by using the Cauchy--Schwarz inequality, the
hypothesis about the sequence $(t_k)_{k}$, and the weak convergence of
$u_{t_k}^\varepsilon$ toward $u^\varepsilon$. The support of $\varphi$ is
compact, so we can apply an integration by part to the integral $\int_\mathbb{R}\varphi(x)\frac{\partial}{\partial x}u_{t_k}^\varepsilon(x)\,dx$.
Hence, we obtain
\begin{eqnarray*}
&&\int_\mathbb{R}\varphi(x) \biggl\{\frac{\varepsilon}{2}
\frac{\partial
}{\partial x}u_{t_k}^\varepsilon(x)+u_{t_k}^\varepsilon(x)
\bigl[V'(x)+F'\ast u_{t_k}^\varepsilon(x)
\bigr] \biggr\}\,dx
\\
&&\qquad=\int_\mathbb{R}\varphi(x) \bigl[V'(x)+F'
\ast u_{t_k}^\varepsilon(x) \bigr]u_{t_k}^\varepsilon(x)\,dx-
\int_\mathbb{R}\frac{\varepsilon}{2}\varphi'(x)u_{t_k}^\varepsilon(x)\,dx
 .
\end{eqnarray*}
The weak convergence of $u_{t_k}^\varepsilon$ to $u^\varepsilon$ implies that
the previous term tends toward $\int_\mathbb{R}\varphi (V'+F'\ast
u^\varepsilon )u^\varepsilon-\int_\mathbb{R}\frac{\varepsilon}{2}\varphi'u^\varepsilon$ when $k$ goes to $\infty$. It has already been proven that
$\int_\mathbb{R}\varphi \{\frac{\varepsilon}{2}\frac{\partial
}{\partial x}u_{t_k}^\varepsilon+u_{t_k}^\varepsilon (V'+F'\ast
u_{t_k}^\varepsilon ) \}$ is collapsing when $k$ goes to $\infty
$. We deduce the following statement:
%e8 #&#
%
\renewcommand{\theequation}{\arabic{section}.\arabic{equation}}
\setcounter{equation}{1}
\begin{equation}
\label{eqfrweaksol} \int_\mathbb{R}\varphi
\bigl(V'+F'\ast u^\varepsilon \bigr)u^\varepsilon-
\int_\mathbb{R}\frac{\varepsilon}{2}\varphi'u^\varepsilon=0
 .
\end{equation}
\item[\textit{Step} 4.] This means that $u^\varepsilon$ is a weak solution of the equation
\[
\frac{\varepsilon}{2}\frac{\partial}{\partial x}u(x)+ \bigl[V'(x)+F'
\ast u(x) \bigr]u(x)=0 .
\]
Now, we consider a smooth function $\widetilde{\varphi}$ with compact
support $[a,b]$. We put
\[
\varphi(x):=\exp \biggl\{\frac{2}{\varepsilon} \bigl[V(x)+F\ast u^\varepsilon
(x) \bigr] \biggr\}\widetilde{\varphi}'(x) .
\]
$\varphi$ is also a smooth function with compact support. Indeed, the
application $x\mapsto F\ast u^\varepsilon(x)$ is a polynomial function
parametrized by the moments of $u^\varepsilon$, and these moments are
bounded. Equality \eqref{eqfrweaksol} becomes
\[
\int_\mathbb{R}\widetilde{\varphi}''(x)
\exp \biggl\{\frac{2}{\varepsilon
} \bigl[V(x)+F\ast u^\varepsilon(x) \bigr] \biggr
\}u^\varepsilon(x)\,dx=0 .
\]
By applying Weyl's lemma, we deduce that $\exp [\frac{2}{\varepsilon
} (V+F\ast u^\varepsilon ) ]u^\varepsilon$ is a smooth
function. Moreover, its second derivative is equal to $0$. Then, there
exists $A,B\in\mathbb{R}$ such that
\[
u^\varepsilon(x)= (Ax+B )\exp \biggl[-\frac{2}{\varepsilon} \bigl(V(x)+F\ast
u^\varepsilon(x) \bigr) \biggr]
\]
for all $x\in\mathbb{R}$. If $A\neq0$, it yields $u^\varepsilon(-Ax)<0$
for $x$ big enough. This is impossible. Consequently, $u^\varepsilon
(x)=Z^{-1} \exp [-\frac{2}{\varepsilon} (V(x)+F\ast u^\varepsilon
(x) ) ]$. This means that $u^\varepsilon$ is a stationary
measure.\quad\qed
\end{longlist}
\noqed\end{pf}
%
%de1.1 #&#
%
\begin{defn}
\label{deffrsetav}
From now on, we call $\mathcal{A}$ the set of the adherence values of
the family $ (u_t^\varepsilon )_{t\in\mathbb{R}_+}$.
\end{defn}
%
%pr1.2 #&#
%
\begin{prop}
\label{propfrlimit}
With the assumptions and the notation of Theorem~\ref{thmfrsubconv},
we have the following limit:
\[
L_0:=\lim_{t\longrightarrow+\infty} \Upsilon_\varepsilon
\bigl(u_t^\varepsilon \bigr)=\Upsilon_\varepsilon
\bigl(u^\varepsilon\bigr) .
\]
\end{prop}
\begin{pf}
The convergence from the quantity $\int_\mathbb{R} Vu_{t_k}^\varepsilon
+\frac{1}{2}\int_\mathbb{R} (F\ast u_{t_k}^\varepsilon
)u_{t_k}^\varepsilon$ toward
$\int_\mathbb{R} Vu^\varepsilon+\frac{1}{2}\int_\mathbb{R} (F\ast u^\varepsilon )u^\varepsilon$ is a consequence of
Theorem~\ref{thmfrsubconv}. So we focus on the entropy term.

First of all, we aim to prove that $ (u_{t_k}^\varepsilon )_k$ is
uniformly bounded in the space $W^{1,1}$. For doing this, we will bound
the integral on $\mathbb{R}$ of $\frac{\partial}{\partial
x}u_{t_k}^\varepsilon(x)$. The triangular inequality provides
\[
\int_\mathbb{R}\biggl\llvert \frac{\partial}{\partial x}u_{t_k}^\varepsilon
\biggr\rrvert \leq\frac{2}{\varepsilon}\int_\mathbb{R}\llvert
\eta_{t_k}\rrvert +\frac
{2}{\varepsilon}\int_\mathbb{R}
\bigl\llvert V'+F'\ast u_{t_k}^\varepsilon
\bigr\rrvert u_{t_k}^\varepsilon,
\]
where $t\mapsto\eta_t$ is defined in Definition~\ref{defintroxi}. By
using \eqref{majoration} and the growth property of $V'$ and $F'$, it yields
\begin{eqnarray*}
\int_\mathbb{R}\bigl\llvert V'(x)+F'
\ast u_{t_k}^\varepsilon(x)\bigr\rrvert u_{t_k}^\varepsilon(x)\,dx
\leq C_1\int_\mathbb{R} \bigl(1+\bigl|x^{2q}\bigr|
\bigr)u_{t_k}^\varepsilon(x)\,dx\leq C_2,
\end{eqnarray*}
where $C_2$ is a constant. By using the Cauchy--Schwarz inequality,
like in the proof of Theorem~\ref{thmfrsubconv}, we obtain
\[
\int_\mathbb{R}\bigl\llvert \eta_{t_k}(x)\bigr\rrvert
\,dx\leq\sqrt{-\xi'(t_k)} .
\]
The quantity $\sqrt{-\xi'(t_k)}$ tends toward $0$, so it is bounded.
Finally, it leads to
\[
\int_\mathbb{R}\biggl\llvert \frac{\partial}{\partial x}u_{t_k}^\varepsilon
(x)\biggr\rrvert \,dx\leq C_3,
\]
where $C_3$ is a constant. Consequently, $u_{t_k}^\varepsilon(x)\leq
u_{t_k}^\varepsilon(0)+C$ for all $x\in\mathbb{R}$. And, since the
sequence $ (u_{t_k}^\varepsilon(0) )_{k}$ converges, it is
bounded, so there exists a constant $C_4$ such that $u_{t_k}^\varepsilon
(x)\leq C_4$ for all $k\in\mathbb{N}$ and $x\in\mathbb{R}$. It is then
easy to prove the convergence of $\int_\mathbb{R} u_{t_k}^\varepsilon
(x)\log (u_{t_k}^\varepsilon(x) )\,dx$ toward $\int_\mathbb{R}
u^\varepsilon(x)\log (u^\varepsilon(x) )\,dx$.

Indeed, the application $x\mapsto u_{t_k}^\varepsilon(x)\log
(u_{t_k}^\varepsilon(x) )$ is lower-bounded, uniformly with respect
to $k$. We can then apply the Lebesgue theorem which provides the
convergence---when $k$ goes to infinity---of $\int_\mathbb{R}
u_{t_k}^\varepsilon(x)\log (u_{t_k}^\varepsilon(x) )\mathds{1}_{\{
|x|\leq R\}}\,dx$ toward $\int_\mathbb{R} u^\varepsilon(x)\log
(u^\varepsilon(x) )\mathds{1}_{\{|x|\leq R\}}\,dx$ for all
$R\geq0$.\vadjust{\goodbreak}
The other integral is split into two terms. The first one is
\begin{eqnarray*}
\int_\mathbb{R} u_{t_k}^\varepsilon(x)\log
\bigl(u_{t_k}^\varepsilon(x) \bigr)\mathds{1}_{\{|x|>R  ;
u_{t_k}^\varepsilon(x)\geq1\}}\,dx&\leq&
\log (C)u_{t_k}^\varepsilon \bigl( [-R;R ]^c \bigr)
\\
&\leq&\frac{\log(C) M_0}{R^2} .
\end{eqnarray*}
The second term is bounded as in the proof of Lemma
\ref{lemfrminoration}:
\begin{eqnarray*}
&&-\int_\mathbb{R} u_{t_k}^\varepsilon(x)\log
\bigl(u_{t_k}^\varepsilon(x) \bigr)\mathds{1}_{\{|x|>R  ;  u_{t_k}^\varepsilon(x)<1\}}\,dx
\\
&&\qquad\leq\int_{[-R;R]^c} \bigl\{|x|u_{t_k}^\varepsilon(x)-
\gamma \bigl(u_{t_k}^\varepsilon(x) \bigr)e^{-{|x|}/{2}} \bigr\}\,dx
\leq\frac
{M_0}{R}+4e^{-{R}/{2}} .
\end{eqnarray*}
Consequently, $\Upsilon_\varepsilon (u_{t_k}^\varepsilon )$
converges toward $\Upsilon_\varepsilon(u^\varepsilon)$, then $\Upsilon_\varepsilon (u_t^\varepsilon )$ converges toward $\Upsilon_\varepsilon
(u^\varepsilon)$ since the free-energy is monotonous.

By taking $R$ big enough and then $k$ big enough, we can make the
following quantity arbitrarily small: $\llvert \int u_{t_k}^\varepsilon\log
(u_{t_k}^\varepsilon )-\int u^\varepsilon\log (u^\varepsilon
)\rrvert $.
\end{pf}
%
%s1.3 #&#
\subsection{Consequences}
When $V$ is symmetric, Proposition 3.1 (resp., Theorem~4.6) in \cite
{HT1} states the existence of at least three stationary measures for
$\varepsilon$ small enough if $F'$ is linear [resp., if
$ \sum_{p=0}^\infty\frac{\llvert F^{(p+2)}(a)\rrvert }{p!}a^p<F''(0)+V''(a)$].
Theorem~\ref{thmfrsubconv} permits to extend these results.
%co1.3 #&#
%
\begin{cor}
\label{corfrnonuni}
For $\varepsilon$ small enough, process \eqref{eqintroinit} admits at
least three stationary measures: one is symmetric ($u^\varepsilon_0$), and
two are asymmetric ($u^\varepsilon_+$ and $u^\varepsilon_-$). Moreover, for
sufficiently small $\varepsilon$,
$\Upsilon_\varepsilon(u^\varepsilon_+)=\Upsilon_\varepsilon(u^\varepsilon_-)<\Upsilon_\varepsilon(u^\varepsilon_0)$.
\end{cor}
\begin{pf}
We know by Theorem 4.5 in~\cite{HT1} that there exists a symmetric
stationary measure $u^\varepsilon_0$. Theorem 5.4 in~\cite{HT2} implies
the weak convergence of $u^\varepsilon_0$ toward $\frac{1}{2} (\delta_{x_0}+\delta_{-x_0} )$ in the small-noise limit where $x_0\in
[0;a[$ is the unique solution of
\[
\cases{ %
\displaystyle V'(x_0)+
\frac{1}{2}F'(2x_0)=0,
\vspace*{2pt}\cr
\displaystyle V''(x_0)+\frac{F''(0)}{2}+
\frac{F''(2x_0)}{2}\geq0 .}
\]
Lemma~\ref{lemaenergylimit} provides
\[
\lim_{\varepsilon\to0}\Upsilon_\varepsilon \bigl(u^\varepsilon_0
\bigr)=V(x_0)+\tfrac{1}{4}F(2x_0) \quad\mbox{and}\quad
 \lim_{\varepsilon\to
0}\Upsilon_\varepsilon \bigl(v^\varepsilon_+
\bigr)=V(a)
\]
with
\[
v^\varepsilon_+(x):=Z^{-1}\exp \biggl[-
\frac{2}{\varepsilon
} \bigl(V(x)+F(x-a) \bigr) \biggr] .
\]
We note that $V(x_0)+\frac{1}{4}F(2x_0)>V(a)$. Consequently, for
$\varepsilon$ small enough, we have $\Upsilon_\varepsilon (v^\varepsilon_+ )<\Upsilon_\varepsilon (u^\varepsilon_0 )$.\vadjust{\goodbreak}

We consider now process \eqref{eqintroinit} starting by
$u_0:=v^\varepsilon_+$. This is possible because the $8q^2$th moment of
$v^\varepsilon_+$ is finite. Theorem~\ref{thmfrsubconv} implies the
existence of a sequence $(t_k)_{k}$ which goes to infinity such that
$u_{t_k}^\varepsilon$ converges weakly toward a stationary measure
$u^\varepsilon$ satisfying $\Upsilon_\varepsilon (u^\varepsilon )\leq
\Upsilon_\varepsilon (u_0 )=\Upsilon_\varepsilon (v^\varepsilon_+ )<\Upsilon_\varepsilon (u^\varepsilon_0 )$. So $u^\varepsilon
\neq u^\varepsilon_0$. We immediately deduce the existence of at least two
stationary measures.

If $V''(0)+F''(0)\neq0$, we know by Theorem 1.6 in~\cite{HT3} that
there exists a unique symmetric stationary measure for $\varepsilon$ small
enough. Hence $u^\varepsilon$ is not symmetric.

Let us assume now that $V''(0)+F''(0)=0$. By \eqref{eqfrminoration},
and by the definition of $\Upsilon_\varepsilon^-(u)$, we have
\begin{eqnarray*}
\Upsilon_\varepsilon(u)&\geq&-\frac{\varepsilon}{4}-4\varepsilon e^{-1}+
\int_\mathbb{R} \biggl\{V(x)-\frac{\varepsilon x^2}{4} \biggr\}u(x)\,dx
\\
&&{}+\frac{1}{2}\iint_{\mathbb{R}^2} F(x-y)u(x)u(y)\,dx\,dy
\end{eqnarray*}
for all $u\in\mathcal{M}_{8q^2}$. Since $F''$ is convex, $x\mapsto
F(x)-\frac{F''(0)}{2}x^2$ is nonnegative. It yields
\begin{eqnarray*}
\Upsilon_\varepsilon(u)&\geq&-\frac{\varepsilon}{4}-4\varepsilon e^{-1}+
\int_\mathbb{R} \biggl\{V(x)+\frac{F''(0)}{2}x^2-
\frac{\varepsilon x^2}{4} \biggr\}u(x)\,dx
\\
&&{}-\frac{F''(0)}{2} \biggl(\int_\mathbb{R} xu(x)\,dx
\biggr)^2 .
\end{eqnarray*}
We deduce the following inequality for all the probability measures
satisfying $\int_\mathbb{R} xu(x)\,dx=0$:
\begin{eqnarray*}
\Upsilon_\varepsilon(u)\geq-\frac{\varepsilon}{4}-4\varepsilon e^{-1}+
\int_\mathbb{R} \biggl\{V(x)+\frac{F''(0)}{2}x^2-
\frac{\varepsilon x^2}{4} \biggr\}u(x)\,dx.
\end{eqnarray*}
In particular, this holds for the symmetric measures. Then, for
$\varepsilon$ small enough, $\Upsilon_\varepsilon(u)>\frac{V(a)}{2}$ for all
the symmetric measures. However, $\Upsilon_\varepsilon(v^\varepsilon_+)<\frac
{V(a)}{2}$ [then $\Upsilon_\varepsilon(u^\varepsilon)<\frac{V(a)}{2}$] for
$\varepsilon$ small enough.

Consequently, the process admits at least one asymmetric stationary
measure that we call $u^\varepsilon_+$. The measure
$u^\varepsilon_-(x):=u^\varepsilon_+(-x)$ is invariant too. By construction of
$u^\varepsilon_+$ and $u^\varepsilon_-$,
$\Upsilon_\varepsilon(u^\varepsilon_+)=\Upsilon_\varepsilon^-(u^\varepsilon_-)<\Upsilon_\varepsilon(u^\varepsilon_0)$.
\end{pf}
%
%re1.4 #&#
%
\begin{rem}
\label{remfrasym}
By a similar method, we could also prove the existence of at least one
stationary measure in the asymmetric-landscape case.
\end{rem}
We know by Theorem 3.2 in~\cite{HT1} that if $V''$ is convex and if
$F'$ is linear, there are exactly three stationary measures for
$\varepsilon$ small enough. We present a more general setting. In view of
the convergence, we will prove that the number of relevant stationary
measures is exactly three even if it is a priori possible to imagine
the existence of at least four such measures.\vadjust{\goodbreak}
%th2 #&#
%
\begin{thmm}
\label{trinity}
We assume $F''(0)+V''(0)>0$. Then, for all $M>0$, there exists $\varepsilon
(M)>0$ such that for all $\varepsilon\leq\varepsilon(M)$, the number of
measures $u$ satisfying the two following conditions is exactly three:
\begin{longlist}[(1)]
\item[(1)] $u$ is a stationary measure for the diffusion \eqref{eqintroinit}.
\item[(2)] $\Upsilon_\varepsilon(u)\leq M$.
\end{longlist}
Moreover, if $\deg(V)=2m>2n=\deg(F)$, diffusion \eqref{eqintroinit}
admits exactly three stationary measures for $\varepsilon$ small enough.
\end{thmm}
\begin{pf}
\textit{Plan}. We will begin to prove the second statement (when \mbox{$m>n$}).
For doing this, we will use Corollary~\ref{corfrnonuni} and the
results in~\cite{HT2,HT3}. Then, we will prove the first statement by
using the second one and a minoration of the free-energy for a sequence
of stationary measures which does not verify~(H).

\begin{longlist}[\textit{Step} 3.2.1.]
\item[\textit{Step} 1.] Corollary~\ref{corfrnonuni} implies the existence of
$\varepsilon_0>0$ such that process \eqref{eqintroinit} admits at least
three stationary measures (one is symmetric, and two are asymmetric) if
$\varepsilon<\varepsilon_0$: $u^\varepsilon_+$, $u^\varepsilon_-$ and $u^\varepsilon_0$.

\item[\textit{Step} 2.] First, we assume that $\deg(V)>\deg(F)$.

\item[\textit{Step} 2.1.] Proposition 3.1 in~\cite{HT2} implies that each family
of stationary measures for the self-stabilizing process \eqref
{eqintroinit} verifies condition (H). It has also been shown that
under (H), we can extract a subsequence which converges weakly from any
family of stationary measures $ (u^\varepsilon )_{\varepsilon>0}$ of
the diffusion \eqref{eqintroinit}.

\item[\textit{Step} 2.2.] Since $F''(0)+V''(0)>0$, there are three possible
limiting values: $\delta_0$, $\delta_a$ and $\delta_{-a}$ according to
Proposition 3.7 and Remark 3.8 in~\cite{HT2}.

\item[\textit{Step} 2.3.] As $F''(0)+V''(0)>0$ and $V''$ and $F''$ are convex,
there is a unique symmetric stationary measure for $\varepsilon$ small
enough by Theorem 1.6 in~\cite{HT3}. Also, Theorem 1.6 in~\cite{HT3}
implies there are exactly two asymmetric stationary measures for
$\varepsilon$ small enough. This achieves the proof of the statement.

\item[\textit{Step} 3.] Now, we will prove the first statement. First, if $m>n$,
by applying the second statement, the result is obvious. We assume now
$m\leq n$. Let $M>0$. All the previous results still hold if we
restrict the study to the families of stationary measures which verify
condition (H). Consequently, it is sufficient to show the following
results in order to achieve the proof of the theorem:
\begin{longlist}[(1)]
\item[(1)] $\sup \{\Upsilon_\varepsilon (u^\varepsilon_0 ) ;
\Upsilon_\varepsilon (u^\varepsilon_+ ) ; \Upsilon_\varepsilon
(u^\varepsilon_- ) \}<M$ for $\varepsilon$ small enough.
\item[(2)] If $ (u^{\varepsilon_k} )_{k}$ is a sequence of
stationary measures, $\int_\mathbb{R} x^{2n}u^{\varepsilon_k}(x)\,dx\to\infty
$ implies $\Upsilon_{\varepsilon_k} (u^{\varepsilon_k} )\to\infty$.
\end{longlist}
\item[\textit{Step} 3.1.] Lemma~\ref{lemaenergylimit} tells us that $\Upsilon_\varepsilon (u^\varepsilon_0 )$ [resp., $\Upsilon_\varepsilon
(u^\varepsilon_+ )=\Upsilon_\varepsilon (u^\varepsilon_- )$] tends
toward $0$ [resp., $V(a)<0$] when $\varepsilon$ goes to $0$. Hence, the
first point is obvious.

\item[\textit{Step} 3.2.] We will prove the second point. We recall lower-bound
\eqref{eqfrminoration},
\[
\Upsilon_\varepsilon^-(u)\geq-\frac{\varepsilon}{4}-4\varepsilon e^{-1}+
\int_\mathbb{R} \biggl(V(x)-\frac{\varepsilon}{4}x^2
\biggr)u(x)\,dx .
\]
As $V(x)\geq C_4x^4-C_2x^2$ and $\Upsilon_\varepsilon^-(u)\leq\Upsilon_\varepsilon(u)$ for all smooth $u$, we obtain
\[
\Upsilon_\varepsilon(u)\geq\int_\mathbb{R}
x^2u(x)\,dx-C,
\]
where $C$ is a constant. It is now sufficient to prove that $\int_\mathbb{R} x^{2n}u^{\varepsilon_k}(x)\,dx\to\infty$ implies $\int_\mathbb
{R} x^{2}u^{\varepsilon_k}(x)\,dx\to\infty$. We will not write the index $k$
for simplifying the reading. We proceed areductio ad absurdum by
assuming the existence of a sequence $ (u^\varepsilon )_\varepsilon$
which verifies $\int_\mathbb{R} x^{2n}u^{\varepsilon}(x)\,dx\to\infty$ and
$\int_\mathbb{R} x^{2}u^{\varepsilon}(x)\,dx\to C_+\in\mathbb{R}_+$.

\item[\textit{Step} 3.2.1.] By taking the notation of~\cite{HT2}, we have the
equality $u^\varepsilon(x)=Z^{-1}\exp [-\frac{2}{\varepsilon}
(W_\varepsilon(x) ) ]$ with
\begin{eqnarray*}
W_\varepsilon(x)&:=&V(x)+F\ast u^\varepsilon(x)=\sum
_{k=1}^{2n}\omega_k(\varepsilon)
 x^k  \mbox{,}
\\
\omega_k(\varepsilon)&:=&\frac{1}{k!} \Biggl\{V^{(k)}(0)+(-1)^k
\sum_{j\geq
{k}/{2}}^{2n}\frac{F^{(2j)}(0)}{(2j-k)!}m_{2j-k}(
\varepsilon) \Biggr\}\quad\mbox{and}
\\
m_{l}(\varepsilon)&:=&\int_\mathbb{R}
x^lu^\varepsilon (x)\,dx\qquad \forall l\in\mathbb{N} .
\end{eqnarray*}
We introduce $\omega(\varepsilon):=\sup \{\llvert \omega_k(\varepsilon
)\rrvert^{{1}/{(2n-k)}}  ;  1\leq k\leq2n \}$.

\item[\textit{Step} 3.2.2.] We note that $\omega_{2n}(\varepsilon)=\frac
{V^{(2n)}(0)+F^{(2n)}(0)}{(2n)!}>0$. Then, $\omega(\varepsilon)$ is
uniformly lower-bounded. Consequently, we can divide by $\omega(\varepsilon
)$.

\item[\textit{Step} 3.2.3.] The change of variable $x:=\omega(\varepsilon)y$ provides
\[
\frac{m_{2l}(\varepsilon)}{\omega(\varepsilon)^{2l}}=\frac{\int_\mathbb{R}
y^{2l}\exp [-({2}/{\widehat{\varepsilon}})\widehat{W}^\varepsilon
(y) ]\,dy}{\int_\mathbb{R}\exp [-{2}/({\widehat{\varepsilon
}})\widehat{W}^\varepsilon(y) ]\,dy} \qquad\mbox{with } \widehat
{W}^\varepsilon(x):=\sum_{k=1}^{2n}
\frac{\omega_k(\varepsilon)}{\omega
(\varepsilon)^{2n-k}}x^k
\]
for all $l\in\mathbb{N}$, with $\widehat{\varepsilon}:=\frac{\varepsilon
}{\omega(\varepsilon)^{2n}}$.

\item[\textit{Step} 3.2.4.] The $2n$ sequences\vspace*{-1pt} $ (\frac{\omega_k(\varepsilon
)}{\omega(\varepsilon)^{2n-k}} )_\varepsilon$ are bounded so we can
extract a subsequence of $\varepsilon$\vspace*{-1pt} (we continue to write $\varepsilon$
for simplifying) such that $\frac{\omega_k(\varepsilon)}{\omega(\varepsilon
)^{2n-k}}$ converges toward $\widehat{\omega}_k$ when $\varepsilon\to0$.
We put $\widehat{W}(x):=\sum_{k=1}^{2n}\widehat{\omega}_kx^k$. We call
$A_1,\ldots,A_r$ the $r\geq1$ location(s) of the global minimum of
$\widehat{W}$.

\item[\textit{Step} 3.2.5.] By applying the result of Lemma~\ref{lemalaplace},
we can extract a subsequence (and we continue to denote it by $\varepsilon
$) such that $\frac{\int_\mathbb{R} y^{2l}\exp [-{2}/{\widehat
{\varepsilon}}\widehat{W}^\varepsilon(y) ]\,dy}{\int_\mathbb{R}\exp
[-{2}/{\widehat{\varepsilon}}\widehat{W}^\varepsilon(y) ]\,dy}$
converges toward $\sum_{j=1}^rp_jA_j^{2l}$ where $p_1+\cdots+p_r=1$ and
$p_j\geq0$.

\item[\textit{Step} 3.2.6.] If $ (\omega(\varepsilon) )_\varepsilon$ is
bounded, since the quantity $\sum_{j=1}^rp_jA_j^{2n}$ is finite, we
deduce that $ (m_{2n}(\varepsilon) )_{\varepsilon}$ is bounded too.
Since $m_{2n}(\varepsilon)$ tends toward infinity when $\varepsilon$ goes to
$0$, we deduce that $(\omega(\varepsilon))_\varepsilon$ converges toward
infinity. As $m_2(\varepsilon)$ is bounded, the quantity $\frac
{m_2(\varepsilon)}{\omega(\varepsilon)^2}$ vanishes when $\varepsilon$ goes to
$0$. This means $\sum_{j=1}^rp_jA_j^2=0$ which implies $A_j=0$ for all
$1\leq j\leq r$. Then $\sum_{j=1}^rp_jA_j^{2n}=0$. Consequently,
$m_{2n}(\varepsilon)=o \{\omega(\varepsilon)^{2n} \}$. The Jensen
inequality provides $m_k(\varepsilon)=o \{\omega(\varepsilon)^{k} \}
$.

\item[\textit{Step} 3.2.7.] We recall the definition of $\omega_k(\varepsilon)$,
\[
\omega_k(\varepsilon)=\frac{1}{k!} \Biggl\{V^{(k)}(0)+(-1)^k
\sum_{j\geq
{k}/{2}}^{2n}\frac{F^{(2j)}(0)}{(2j-k)!}m_{2j-k}(
\varepsilon) \Biggr\} .
\]
We deduce $\omega_k(\varepsilon)=O \{m_{2n-k}(\varepsilon) \}=o
\{\omega(\varepsilon)^{2n-k} \}$. So
\[
\omega(\varepsilon)=\sup \bigl\{\bigl\llvert \omega_k(\varepsilon)\bigr
\rrvert^{
{1}/{(2n-k)}}  ;  1\leq k\leq2n \bigr\}=o \bigl\{\omega(\varepsilon)
\bigr\} .
\]
This is a contradiction which achieves the proof.\quad\qed
\end{longlist}
\noqed\end{pf}
This theorem means that---even if the diffusion \eqref{eqintroinit}
admits more than three stationary measures---there are only three
stationary measures which play a role in the convergence. Indeed, if we
take a measure $u_0$ with a finite free-energy, we know that for
$\varepsilon$ small enough, there are only three (maybe fewer) stationary
measures which can be adherence value of the family $ (u_t^\varepsilon
)_{t\in\mathbb{R}_+}$.

The assumption (LIN) implies (M3) (and (M3)$'$ because it is weaker) and
(0M1) for $\varepsilon$ small enough. The condition (SYN) implies (M3)$'$
and (0M1) for $\varepsilon$ small enough. Furthermore, if $\deg(V)>\deg
(F)$, (SYN) implies (M3) when $\varepsilon$ is less than a threshold.

This description of the stationary measures permits us to obtain the
principal result, that is to say, the long-time convergence of the process.
%s2 #&#
\section{Global convergence}
%s2.1 #&#
\subsection{Statement of the theorem}
We write the main result of the paper:
%th3 #&#
%
\begin{thmm}
\label{convergence}
Let $du_0$ be a probability measure which verifies \textup{(FE)} and~\textup{(FM)}. Under
\textup{(M3)}, $u_t^\varepsilon$ converges weakly toward a stationary measure.
\end{thmm}
The proof is postponed in Section~\ref{subsectionproof}. First, we
will discuss briefly the assumptions.

%s2.2 #&#
\subsection{Remarks on the assumptions}
\subsubsection*{$du_0$ is absolutely continuous with respect to the
Lebesgue measure}
We shall use Theorem~\ref{thmfrsubconv} and prove that the family
$ (u_t^\varepsilon )_{t\in\mathbb{R}_+}$ admits a unique
adherence value. This theorem requires that the initial law is
absolutely continuous with respect to the Lebesgue measure.\vadjust{\goodbreak} However, it
is possible to relax this hypothesis by using the following result (see
Lemma 2.1 in~\cite{HT1} for a proof):

\textit{Let $du_0$ be a probability measure which verifies $\int_\mathbb
{R} x^{8q^2}\,du_0(x)<+\infty$. Then, for all $t>0$, the probability
$du_t^\varepsilon$ is absolutely continuous with respect to the Lebesgue
measure.}

Consequently, it is sufficient to apply Theorem~\ref{convergence} to
the probability measure~$u_1^\varepsilon$ since there is a unique solution
to the nonlinear equation \eqref{eqintroinit}.

\subsubsection*{The entropy of $du_0$ is finite}
An essential point of the proof is the convergence of the free-energy.
To be sure of this, we assume that it is finite at time~$0$. The
assumption about the moments implies $\Upsilon_\varepsilon
(u_t^\varepsilon )<+\infty$ if and only if $\int_\mathbb{R}
u_t^\varepsilon\log (u_t^\varepsilon )<+\infty$.

If $V$ was convex, a little adaptation of the theorem in~\cite{OV2001}
(taking into account the fact that the drift is not homogeneous here)
would provide the nonoptimal following inequality:
\[
\Upsilon_\varepsilon \bigl(u_t^\varepsilon \bigr)\leq
\frac{1}{2t}\inf \bigl\{ \sqrt{\mathbb{E}\llvert X-Y\rrvert^2}
 ;  \mathcal{L} (X )=u_t^\varepsilon  ;  \mathcal{L} (Y
)=v_t^\varepsilon \bigr\}
\]
with
\[
v_t^\varepsilon(x):=Z^{-1}\exp
\biggl[-\frac{2}{\varepsilon
} \bigl(V(x)+F\ast u_t^\varepsilon(x)
\bigr) \biggr]
\]
for all $t>0$. The second moment of $u_t^\varepsilon$ is upper-bounded
uniformly with respect to $t$. By using the convexity of $V$ and $F$,
we can prove the same thing for $v_t^\varepsilon$. Consequently, since
$t>0$, the free-energy is finite so the entropy is finite. However, in
this paper, we deal with nonconvex landscape, so we will not relax
this hypothesis.

\subsubsection*{All the moments are finite}
Theorem~\ref{thmfrsubconv} tells us that we can extract a sequence
from the family $ (u_t^\varepsilon )_{t\in\mathbb{R}_+}$ such
that it converges toward a stationary measure. The last step in order
to obtain the convergence is the uniqueness of the limiting value. The
most difficult part will be to prove this uniqueness when the symmetric
stationary measure $u^\varepsilon_0$ is an adherence value and the only
one of these adherence values to be stationary. To do this, we will
consider a function like this one:
\[
\Phi(u):=\int_\mathbb{R} \varphi(x)u(x)\,dx,
\]
where $\varphi$ is an odd and smooth function with compact support such
that $\varphi(x)=x^{2l+1}$ for all $x$ in a compact subset of $\mathbb
{R}$. Then, we will prove---by proceeding a \textit{reductio ad
absurdum}---that there exists an integer $l$ such that $\Phi
(u^\varepsilon_0 )\neq\Phi (u_\infty^\varepsilon )$, where
$u_\infty^\varepsilon$ would be another limiting value which is a
stationary measure. This inequality will allow us to construct a
stationary measure $u^\varepsilon$ such that $\Phi(u^\varepsilon)\notin \{
\Phi (u^\varepsilon_0 ) ; \Phi (u^\varepsilon_+ ) ; \Phi
(u^\varepsilon_- ) \}$. This implies the existence of a
stationary measure which does not belong to $ \{u^\varepsilon_0 ;
u^\varepsilon_+ ; u^\varepsilon_- \}$. Under (M3), it is impossible.\vadjust{\goodbreak}

We make the integration with an ``almost-polynomial'' function because
we need the square of the derivative of such function to be uniformly
bounded with respect to the time.

However, it is possible to relax the condition (FM). Indeed, according
to Proposition~\ref{propafm}, if we assume that $\int_\mathbb{R}
x^{8q^2}\,du_0(x)<+\infty$ (the condition used for the existence of a
strong solution), we have
\[
\int_\mathbb{R} x^{2l}u_t^\varepsilon(x)\,dx<+
\infty \qquad\forall t>0, l\in \mathbb{N} .
\]

\subsubsection*{Hypothesis \textup{(M3)}}
As written before, the key for proving the uniqueness of the adherence
value is to proceed a reductio ad absurdum and then to construct a
stationary measure $u^\varepsilon$ such that $\Phi(u^\varepsilon)$ takes a
forbidden value [a value different from $\Phi (u^\varepsilon_0
)$, $\Phi (u^\varepsilon_+ )$ and $\Phi (u^\varepsilon_-
)$].

But, it is possible to deal with a weaker hypothesis. Indeed, by
considering an initial law with finite free-energy and since the
free-energy is decreasing, it is impossible for $u_t^\varepsilon$ to
converge toward a stationary measure with a higher energy.
Consequently, we can consider (M3)$'$ instead of (M3).

All of these remarks allow us to obtain the following result:
%th4 #&#
%
\begin{thmm}
\label{convergence2}
Let $du_0$ be a probability measure with finite entropy. If $V$ and~$F$
are polynomial functions such that $F''(0)+V''(0)>0$, $u_t^\varepsilon$
converges weakly toward a stationary measure for $\varepsilon$ small enough.
\end{thmm}
%
%s2.3 #&#
\subsection{Proof of the theorem}
\label{subsectionproof}
In order to obtain the statement of Theorem~\ref{convergence}, we will
provide two lemmas and one proposition about the free-energy. The
lemmas state that a probability measure which verifies simple
properties and with a level of energy is necessary a stationary measure
for the self-stabilizing process \eqref{eqintroinit}. The third one
allows us to confine all the adherence values under a level of energy.
%le2.1 #&#
%
\begin{lem}
\label{lemgcminimizer}
Under \textup{(M3)}, if $u$ is a probability measure which satisfies \textup{(FE)} and
\textup{(ES)}, the inequality $\Upsilon_\varepsilon(u)\leq\Upsilon_\varepsilon
(u^\varepsilon_\pm)$ implies $u\in \{u^\varepsilon_+ ; u^\varepsilon_- \}$.
\end{lem}
\begin{pf}
Let $u$ be such a measure. We consider the process \eqref
{eqintroinit} starting by the initial law $u_0:=u$. Theorem \ref
{thmfrsubconv} implies that there exists a stationary measure
$u^\varepsilon$ such that $\Upsilon_\varepsilon (u_t^\varepsilon )$
converges toward $\Upsilon_\varepsilon(u^\varepsilon)$.

However, according to Propositions~\ref{propintrodecrease} and~\ref{propfrlimit},
\[
\Upsilon_\varepsilon\bigl(u^\varepsilon\bigr)=\lim_{t\to+\infty}
\Upsilon_\varepsilon \bigl(u_t^\varepsilon \bigr)\leq
\Upsilon_\varepsilon \bigl(u_t^\varepsilon \bigr)\leq
\Upsilon_\varepsilon(u)\leq\Upsilon_\varepsilon\bigl(u^\varepsilon_\pm
\bigr) .
\]
Condition (M3) provides $u^\varepsilon\in \{u^\varepsilon_+;  u^\varepsilon_- ; u^\varepsilon_0 \}$.
But, $\Upsilon_\varepsilon(u^\varepsilon)\leq
\Upsilon_\varepsilon(u^\varepsilon_\pm)<\Upsilon_\varepsilon(u^\varepsilon_0)$ so
$u^\varepsilon\in \{u^\varepsilon_+ ; u^\varepsilon_- \}$. Without
loss of generality, we will assume $u^\varepsilon=u^\varepsilon_+$.

Consequently, the function $\xi$ (see Definition~\ref{defintroxi}) is
constant. We deduce that $\xi'(t)=0$ for all $t\geq0$. Lemma \ref
{lemfrstation} implies that $u_t^\varepsilon$ is a stationary measure.
This means that $u=u_0=u^\varepsilon=u^\varepsilon_+$.\vadjust{\goodbreak}
\end{pf}
We have a similar result with the symmetric measures:
%le2.2 #&#
%
\begin{lem}
\label{lemgcminimizersym}
Under \textup{(0M1)}, if $u$ is a symmetric probability measure satisfying~\textup{(FE)}
and \textup{(ES)}, $\Upsilon_\varepsilon(u)\leq\Upsilon_\varepsilon(u^\varepsilon_0)$
implies $u=u^\varepsilon_0$.
\end{lem}
The key-argument is the following: if the initial law is symmetric,
then the law at time $t$ is still symmetric. The proof is similar to
the previous one, so it is left to the reader's attention.

Before making the convergence, we need a last result on the adherence
values: \textit{the free-energy of a limiting value is less than the
limit value of the free-energy.}
%pr2.3 #&#
%
\begin{prop}
\label{propgcfatou}
We assume that $u_\infty^\varepsilon$ is an adherence value of the family
$ (u_t^\varepsilon )_{t\in\mathbb{R}_+}$. We call
$L_0:= \lim_{t\to+\infty}\Upsilon_\varepsilon (u_t^\varepsilon
)$. Then $\Upsilon_\varepsilon(u_\infty^\varepsilon)\leq L_0$.
\end{prop}
\begin{pf}
As $u_\infty^\varepsilon$ is an adherence value of the family $
(u_t^\varepsilon )_{t\in\mathbb{R}_+}$, there exists an increasing
sequence $(t_k)_{k}$ which goes to infinity such that $u_{t_k}^\varepsilon
$ converges weakly toward $u_\infty^\varepsilon$. We remark
\begin{eqnarray*}
\Upsilon \bigl(u_{t_{k}}^\varepsilon \bigr)&=&V(a)+\int
_\mathbb{R} \bigl(V(x)-V(a) \bigr)u_{t_{k}}^\varepsilon(x)\,dx
\\
&&{}+\frac{1}{2}\iint_{\mathbb{R}^2}F(x-y)u_{t_{k}}^\varepsilon
(x)u_{t_{k}}^\varepsilon(y)\,dx\,dy,
\end{eqnarray*}
where the functional $\Upsilon$ is defined in \eqref
{eqintrometapotentiel}. As $V(x)-V(a)\geq0$ for all $x\in\mathbb{R}$,
the Fatou lemma implies $ \Upsilon (u_\infty^\varepsilon
)\leq\liminf_{k\to\infty}\Upsilon (u_{t_{k}}^\varepsilon
)$. In the same way,
\begin{eqnarray*}
&&\int_\mathbb{R} u_\infty^\varepsilon(x)\log
\bigl(u_\infty^\varepsilon (x) \bigr)\mathds{1}_{\{u_\infty^\varepsilon(x)\geq1\}}\,dx
\\
&&\qquad\leq\liminf_{k\to\infty}\int_\mathbb{R}
u_{t_{k}}^\varepsilon(x)\log \bigl(u_{t_{k}}^\varepsilon(x)
\bigr)\mathds{1}_{\{u_{t_{k}}^\varepsilon(x)\geq1\}
}\,dx .
\end{eqnarray*}
Let $R>0$. By putting $\gamma_k^-(x):=u_{t_{k}}^\varepsilon(x)\log
(u_{t_{k}}^\varepsilon(x) )\mathds{1}_{\{u_{t_{k}}^\varepsilon(x)<1\}
}\mathds{1}_{\{|x|\leq R\}}$, we note that $\llvert \gamma_k^-(x)\rrvert \leq e^{-1}\mathds{1}_{\{|x|\leq R\}}$
for all $x\in\mathbb{R}$ and
$k\in\mathbb{N}$. We can apply the Lebesgue theorem,
\begin{eqnarray*}
\int_\mathbb{R} u_\infty^\varepsilon(x)\log
\bigl(u_\infty^\varepsilon(x) \bigr)\mathds{1}_{\{u_\infty^\varepsilon(x)\geq1\}}
\mathds{1}_{\{|x|\leq R\}
}=\lim_{k\to\infty}\int_\mathbb{R}
\gamma_k^-(x)\,dx .
\end{eqnarray*}
We put $\gamma_k^+(x):=u_{t_{k}}^\varepsilon(x)\log (u_{t_{k}}^\varepsilon
(x) )\mathds{1}_{\{u_{t_{k}}^\varepsilon(x)<1\}}\mathds{1}_{\{|x|>R\}
}$. By proceeding as in the proof of Lemma~\ref{lemfrminoration}, we have
\begin{eqnarray*}
-\gamma_k^+(x)&=&-u_{t_{k}}^\varepsilon(x)\log
\bigl(u_{t_{k}}^\varepsilon (x) \bigr)\mathds{1}_{\{e^{-|x|}\leq u_{t_{k}}^\varepsilon(x)<1\}}\mathds
{1}_{\{|x|>R\}}
\\
&&{}-u_{t_{k}}^\varepsilon(x)\log \bigl(u_{t_{k}}^\varepsilon(x)
\bigr)\mathds {1}_{\{u_{t_{k}}^\varepsilon(x)<e^{-|x|}\}}\mathds{1}_{\{|x|>R\}}
\\
&\leq&|x|u_{t_{k}}^\varepsilon(x)\mathds{1}_{\{|x|>R\}}+2e^{-1}e^{-{|x|}/{2}}
\mathds{1}_{\{|x|>R\}} .
\end{eqnarray*}
Consequently, it leads to the lower-bound
\begin{eqnarray*}
\int_\mathbb{R} u_{t_{k}}^\varepsilon(x)\log
\bigl(u_{t_{k}}^\varepsilon (x) \bigr)\mathds{1}_{\{u_{t_{k}}^\varepsilon(x)<1\}}
\mathds{1}_{\{|x|>R\}
}\,dx\geq-\frac{M_0}{R}-8e^{-1}e^{-{R}/{2}},
\end{eqnarray*}
where $M_0$ is defined in \eqref{majoration}.

By introducing $\widehat{\Upsilon}_\varepsilon(u):=\Upsilon_\varepsilon
(u)-\frac{\varepsilon}{2}\int_\mathbb{R} u(x)\log(u(x))\mathds{1}_{\{
u(x)<1\}}\mathds{1}_{\{|x|>R\}}\,dx$, we obtain
\begin{eqnarray*}
\Upsilon_\varepsilon \bigl(u_\infty^\varepsilon \bigr)\leq
\widehat{\Upsilon}_\varepsilon \bigl(u_\infty^\varepsilon \bigr)&
\leq&\liminf_{k\to\infty
}\widehat{\Upsilon}_\varepsilon \bigl(u_{t_k}^\varepsilon
\bigr)
\\
&\leq&\liminf_{k\to\infty}\Upsilon_\varepsilon \bigl(u_{t_k}^\varepsilon
\bigr)+\frac{M_0\varepsilon}{2R}+4e^{-1}\exp \biggl(-\frac{R}{2}
\biggr)\varepsilon
\\
&\leq& L_0+\frac{M_0\varepsilon}{2R}+4e^{-1}\varepsilon\exp \biggl(-
\frac
{R}{2} \biggr)
\end{eqnarray*}
for all $R>0$. Consequently, $\Upsilon_\varepsilon(u_\infty^\varepsilon)\leq L_0$.
\end{pf}
\begin{pf*}{Proof of the theorem}
\textit{Plan}: The first step of the proof consists of the application of
the Prohorov theorem since the family of measure is tight. We shall
prove the uniqueness of the adherence value. We will proceed a reductio
ad absurdum. The previous results provide $\mathcal{A}\cap \{
u^\varepsilon_0 ; u^\varepsilon_+ ; u^\varepsilon_- \}\neq\varnothing$
where $\mathcal{A}$ is introduced in Definition~\ref{deffrsetav}. We
will then study all the possible cases, and we will prove that all of
these cases imply contradictions. If $\mathcal{A}\cap \{
u^\varepsilon_0 ; u^\varepsilon_+ ; u^\varepsilon_- \}= \{u^\varepsilon_+ \}$ and $\mathcal{A}\cap \{u^\varepsilon_0 ; u^\varepsilon_+ ; u^\varepsilon_- \}= \{u^\varepsilon_- \}$ imply
contradiction since $u^\varepsilon_+$ and $u^\varepsilon_-$ are the unique
minimizers of the free-energy. The cases $u^\varepsilon_0\in\mathcal{A}$
and $\mathcal{A}\cap \{u^\varepsilon_0 ; u^\varepsilon_+ ;
u^\varepsilon_- \}= \{u^\varepsilon_+ ; u^\varepsilon_- \}$
contradict (M3).
\begin{longlist}[\textit{Step} 2.1.1.]
\item[\textit{Step} 1.] Inequality \eqref{majoration} implies that the family of
probability measures $ \{u_t^\varepsilon ; t\in\mathbb{R}_+ \}$
is tight. Prohorov's theorem allows us to conclude that each extracted
sequence of this family is relatively compact with respect to the weak
convergence. So, in order to prove the statement of the theorem, it is
sufficient to prove that this family admits exactly one adherence
value. We proceed a reductio ad absurdum. We assume in the following
that the family admits at least two adherence values.

\item[\textit{Step} 2.] As condition (M3) is true, there are exactly three
stationary measures: $u^\varepsilon_0$, $u^\varepsilon_+$ and $u^\varepsilon_-$.
By Theorem~\ref{thmfrsubconv}, we know that $\mathcal{A}\cap \{
u^\varepsilon_0 ; u^\varepsilon_+ ; u^\varepsilon_- \}\neq\varnothing$.
We split this step into three cases:
\begin{itemize}
\item$u^\varepsilon_0\in\mathcal{A}$.
\item$\mathcal{A}\cap \{u^\varepsilon_0 ; u^\varepsilon_+ ;
u^\varepsilon_- \}= \{u^\varepsilon_+ ; u^\varepsilon_- \}$.
\item$\mathcal{A}\cap \{u^\varepsilon_0 ; u^\varepsilon_+ ;
u^\varepsilon_- \}= \{u^\varepsilon_+ \}$.
\end{itemize}
By symmetry, we will not deal with the case $\mathcal{A}\cap \{
u^\varepsilon_0 ; u^\varepsilon_+ ; u^\varepsilon_- \}= \{u^\varepsilon_- \}$.

\item[\textit{Step} 2.1.] We will prove that the first case, $u^\varepsilon_0\in
\mathcal{A}$, is impossible. It will be the core of the proof.

\item[\textit{Step} 2.1.1.] Let $u_\infty^\varepsilon$ be an other adherence value
of the family
$ (u_t^\varepsilon )_{t\in\mathbb{R}_+}$. Proposition \ref
{propgcfatou} tells us
$\Upsilon_\varepsilon (u_\infty^\varepsilon )\leq\Upsilon_\varepsilon
(u^\varepsilon_0 )$.
Since $u_\infty^\varepsilon\neq u^\varepsilon_0$, Lemma \ref
{lemgcminimizersym} implies that the law $u_\infty^\varepsilon$ is not
symmetric. We deduce that there exists $l\in\mathbb{N}$ such that $\int_\mathbb{R} x^{2l+1}u_\infty^\varepsilon(x)\,dx\neq0$. Let $R>0$. We
introduce the function
\begin{eqnarray*}
\varphi(x)&:=&x^{2l+1}\mathds{1}_{[-R;R]}(x)
\\
&&{}+x^{2l+1} \mathds{1}_{[R;R+1]}(x) Z^{-1}\int
_x^{R+1}\exp \biggl[-\frac
{1}{(y-R)^2}-
\frac{1}{(y-R-1)^2} \biggr]\,dy
\\
&&{}+x^{2l+1} \mathds{1}_{[-R-1;-R]}(x) Z^{-1}\int
_{-R-1}^x\exp \biggl[-\frac{1}{(y+R)^2}-
\frac{1}{(y+R+1)^2} \biggr]\,dy
\end{eqnarray*}
with
\[
Z:=\int_0^1\exp \biggl[-
\frac{1}{z^2}-\frac{1}{(z-1)^2} \biggr]\,dz .
\]
By construction, $\varphi$ is an odd function, so
$\int_\mathbb{R}\varphi(x)u^\varepsilon_0(x)\,dx=0$. Furthermore,
$\llvert \varphi(x)\rrvert \leq|x|^{2l+1}$. By using the triangular
inequality and (FM), we have
\begin{eqnarray*}
\biggl\llvert \int_\mathbb{R}\varphi(x)u_\infty^\varepsilon(x)\,dx
\biggr\rrvert &\geq&\biggl\llvert \int_\mathbb{R}
x^{2l+1}u_\infty^\varepsilon(x)\,dx\biggr\rrvert -\int
_{[-R;R]^c}\llvert x\rrvert^{2l+1}u_\infty^\varepsilon(x)\,dx
\\
&\geq&\biggl\llvert \int_\mathbb{R} x^{2l+1}u_\infty^\varepsilon(x)\,dx
\biggr\rrvert -\frac
{1}{R^3}C_0,
\end{eqnarray*}
where $C_0:=\sup_{t\in\mathbb{R}_+}\int_\mathbb{R}|x|^{2l+4}u_t^\varepsilon
(x)\,dx<+\infty$.
Since $\int_\mathbb{R} x^{2l+1}u_\infty^\varepsilon(x)\,dx\neq0$, we deduce that
$\int_\mathbb{R}\varphi(x)u_\infty^\varepsilon(x)\,dx\neq0$ for $R$ big
enough. Consequently,
we obtain the existence of a smooth function $\varphi$ with compact
support such that
\[
0=\int_\mathbb{R}\varphi(x)u^\varepsilon_0(x)\,dx<
\int_\mathbb{R}\varphi (x)u_\infty^\varepsilon(x)\,dx
 .
\]
Moreover, we can verify that $\varphi'(x)^2\leq
C(R)x^{4l+2}$ for all $x\in\mathbb{R}$. This implies
$\sup_{t\in\mathbb{R}_+}\int_\mathbb{R}\varphi'(x)^2u_t^\varepsilon
(x)\,dx<+\infty$.

\item[\textit{Step} 2.1.2.] Let $\kappa>0$ such that $\llvert \int_\mathbb{R}\varphi
(x)u^\varepsilon_+(x)\,dx\rrvert >3\kappa$. By definition of~$\mathcal{A}$,
there exist two increasing sequences $(t_k^{(1)})_{k}$ [resp.,
$(t_k^{(2)})_{k}$] which go to infinity such that
$u_{t_k^{(1)}}^\varepsilon$ [resp., $u_{t_k^{(2)}}^\varepsilon$] converges
weakly toward $u^\varepsilon_0$ (resp., $u_\infty^\varepsilon$). We deduce
that there exist two increasing sequences $(r_k)_{k}$ and $(s_k)_{k}$
such that $\int_\mathbb{R}\varphi(x)u_{r_k}^\varepsilon(x)\,dx=\kappa$ and
$\int_\mathbb{R}\varphi(x)u_{s_k}^\varepsilon(x)\,dx=2\kappa$. Then, for all
$k\in\mathbb{N}$, we put $\widehat{r}_k:=\sup \{t\in
[0;s_k ]  \llvert   \int_\mathbb{R}\varphi(x)u_t^\varepsilon
(x)\,dx=\kappa \}$, and then we define $\widehat{s}_k:=\inf \{
s\in [\widehat{r}_k;s_k ]  \llvert  \int_\mathbb
{R}\varphi(x)u_s^\varepsilon(x)\,dx=2\kappa \}$. For simplifying, we
write $r_k$ (resp., $s_k$) instead of $\widehat{r}_k$ (resp., $\widehat
{s}_k$). And we have
\begin{eqnarray*}
\kappa=\int_\mathbb{R}\varphi(x)u_{r_k}^\varepsilon(x)\,dx
\leq\int_\mathbb {R}\varphi(x)u_t^\varepsilon(x)\,dx
\leq\int_\mathbb{R}\varphi (x)u_{s_k}^\varepsilon(x)\,dx=2
\kappa
\end{eqnarray*}
for all $t\in[r_k;s_k]$.\vadjust{\goodbreak}

\item[\textit{Step} 2.1.3.] By applying Proposition~\ref{propaothermeasure}, we
deduce that there exists an increasing sequence $(q_k)_{k}$ going to
$+\infty$ such that $ (u_{q_k}^\varepsilon )_k$ converges weakly
toward a stationary measure $u^\varepsilon$ verifying $\int_\mathbb
{R}\varphi(x)u^\varepsilon(x)\,dx\in[\kappa;2\kappa]$. Since we have the
inequality $\llvert \int_\mathbb{R}\varphi(x)u^\varepsilon_+(x)\,dx\rrvert =\llvert \int_\mathbb{R}\varphi(x)u^\varepsilon_-(x)\,dx\rrvert >3\kappa$, we
deduce $u^\varepsilon=u^\varepsilon_0$. This is impossible since $\int_\mathbb
{R}\varphi(x)u^\varepsilon_0(x)\,dx=0\notin[\kappa;2\kappa]$.

\item[\textit{Step} 2.2.] We deal now with the third case, $\mathcal{A}\cap \{u^\varepsilon_0 ; u^\varepsilon_+ ; u^\varepsilon_- \}= \{
u^\varepsilon_+ ; u^\varepsilon_- \}$.

\item[\textit{Step} 2.2.1.] By definition of $u^\varepsilon_+$ and $u^\varepsilon_-$,
these measures are not symmetric. Consequently, there exists $l\in
\mathbb{N}$ such that $\int_\mathbb{R} x^{2l+1}u^\varepsilon_+(x)\,dx\neq0$.
As $u^\varepsilon_-(x)=u^\varepsilon_+(-x)$, by proceeding as in \textit{Step}~2.1,
we deduce that there exists an increasing sequence $(q_k)_{k\in
\mathbb{N}}$ which goes to $\infty$ such that $u_{q_k}^\varepsilon$
converges weakly toward a stationary measure $u^\varepsilon$ which
verifies $\int_\mathbb{R}\varphi u^\varepsilon\in[\kappa;2\kappa]$ where
$\varphi$ is a smooth function with compact support such that $\int_\mathbb{R}\varphi u^\varepsilon_\pm\notin[\kappa;2\kappa]$. We deduce
that $u^\varepsilon=u^\varepsilon_0$ which contradicts $u^\varepsilon_0\notin
\mathcal{A}$.

\item[\textit{Step} 2.3.] We consider now the last case, $\mathcal{A}\cap \{u^\varepsilon_0 ; u^\varepsilon_+ ; u^\varepsilon_- \}= \{
u^\varepsilon_+ \}$. Proposition~\ref{propfrlimit} implies that
$\Upsilon_\varepsilon (u_t^\varepsilon )$ converges toward $\Upsilon_\varepsilon(u^\varepsilon_+)$. Let $u_\infty^\varepsilon$ be a limit value of
the family $ (u_t^\varepsilon )_{t\in\mathbb{R}_+}$ which is not
$u^\varepsilon_+$. By Proposition~\ref{propgcfatou}, we know that
$\Upsilon_\varepsilon (u_\infty^\varepsilon )\leq\Upsilon_\varepsilon
(u^\varepsilon_+ )= \lim_{t\longrightarrow+\infty
}\Upsilon_\varepsilon (u_t^\varepsilon )$. Then, Lemma \ref
{lemgcminimizer} implies $u_\infty^\varepsilon=u^\varepsilon_-\notin\mathcal
{A}$.\vspace*{-2pt}

\begin{Conclusion*}
The family $ (u_t^\varepsilon )_{t\in\mathbb
{R}_+}$ admits only one adherence value with respect to the weak
convergence. So $u_t^\varepsilon$ converges weakly toward a stationary
measure which achieves the proof.\quad\qed\vspace*{-2pt}
\end{Conclusion*}
\end{longlist}
\noqed\end{pf*}

%s3 #&#
\section{Basins of attraction}
Now we shall provide some conditions in order to precise the limit.\vspace*{-2pt}
%s3.1 #&#
\subsection{\texorpdfstring{Domain of $u^\varepsilon_0$}{Domain of u epsilon 0}}\vspace*{-2pt}
%th5 #&#
%
\begin{thmm}
\label{thmbasym}
Let $du_0$ be a symmetric probability measure which verifies \textup{(FE)} and
\textup{(ES)}. We assume that $V''(0)+F''(0)\neq0$. Then, for $\varepsilon$ small
enough $u_t^\varepsilon$ converges weakly toward $u^\varepsilon_0$.\vspace*{-2pt}
\end{thmm}
\begin{pf}
$V''(0)+F''(0)\neq0$, and both functions $V''$ and $F''$ are convex.
Theorem 1.6 in~\cite{HT3} implies the existence and the uniqueness of a
symmetric stationary measure $u^\varepsilon_0$ for $\varepsilon$ small
enough.

Theorem~\ref{thmfrsubconv} provides the existence of a stationary
measure $u^\varepsilon$ and an increasing sequence $(t_k)_{k}$ which goes
to $\infty$ such that $u_{t_k}^\varepsilon$ converges weakly toward
$u^\varepsilon$ and $\Upsilon_\varepsilon (u_t^\varepsilon )$ converges
toward $\Upsilon_\varepsilon(u^\varepsilon)$. As $u_t^\varepsilon$ is symmetric
for all $t\geq0$, we deduce $u^\varepsilon=u^\varepsilon_0$, the unique
symmetric stationary measure.

We proceed a reductio ad absurdum by assuming the existence of another
sequence $(s_k)_{k}$ which goes to $\infty$ such that $u_{s_k}^\varepsilon
$ does not converge toward $u^\varepsilon_0$. The uniform boundedness of
the second moment with respect to the time permits to extract a
subsequence [we continue to write $(s_k)_{k}$ for simplifying] such
that $u_{s_k}^\varepsilon$ converges weakly toward $u_\infty^\varepsilon\neq
u^\varepsilon_0$. Proposition~\ref{propgcfatou} implies $\Upsilon_\varepsilon(u_\infty^\varepsilon)\leq\Upsilon_\varepsilon(u^\varepsilon_0)$. Lemma
\ref{lemgcminimizersym} implies $u_\infty^\varepsilon=u^\varepsilon_0$.
This is absurd.\vadjust{\goodbreak}
\end{pf}
%
%re3.1 #&#
%
\begin{rem}
\label{remgcboundary}
We assume $V''(0)+F''(0)\neq0$ in order to have a unique symmetric
stationary measure for $\varepsilon$ small enough, that is to say (0M1).
We can extend to the case $V''(0)+F''(0)=0$ by using auniform
propagation of chaos; see Theorem 6.5 in~\cite{TT}. We can also assume
that $n=2$ which means $\deg(F)=4$ by Section 4.2 in~\cite{HT1}.
\end{rem}
%
%re3.2 #&#
%
\begin{rem}
\label{remgcfm}
In the previous theorem, if we assumed (FM) instead of (ES), we could
have directly applied Theorem~\ref{convergence}.
\end{rem}
%
%s3.2 #&#
\subsection{\texorpdfstring{Domain of $u^\varepsilon_\pm$}{Domain of u epsilon +-}}
The principal tool of the previous theorem is the stability of a subset
(all the symmetric measures with a finite $8q^2$-moment). If we could
find an invariant subset which contains $u^\varepsilon_+$, but neither
$u^\varepsilon_0$ nor $u^\varepsilon_-$, we could apply the same method than
previously.

Instead of this, we will consider an inequality linked to the
free-energy and we will exhibit a simple subset included in the domain
of attraction of~$u^\varepsilon_+$. Let us first introduce the following hyperplan:
\begin{eqnarray*}
\mathcal{H}:= \biggl\{u\in\mathcal{C}^\infty (\mathbb{R} ; \mathbb
{R}_+ )   \Big|  \int_\mathbb{R} x^{8q^2}u(x)\,dx<\infty
  \mbox{ and }  \int_\mathbb{R} xu(x)\,dx=0 \biggr\} .
\end{eqnarray*}

%th6 #&#
%
\begin{thmm}
\label{thmbaasym}
Let $du_0$ be a probability measure which verifies \textup{(FE)} and \textup{(FM)}. We
assume also
\[
\Upsilon_\varepsilon (u_0 )<\inf_{u\in\mathcal{H}}
\Upsilon_\varepsilon(u) \quad\mbox{and}\quad \int_\mathbb{R}
xu_0(x)\,dx>0 .
\]
Under \textup{(M3)}, $u_t^\varepsilon$ converges weakly toward $u^\varepsilon_+$.
\end{thmm}
\begin{pf}
We know by Theorem~\ref{convergence} that there exists a stationary
measure $u^\varepsilon$ such that $ (u_t^\varepsilon )_t$ converges
weakly toward $u^\varepsilon$. And, by Proposition~\ref{propfrlimit},
$\Upsilon_\varepsilon (u_t^\varepsilon )$ converges toward $\Upsilon_\varepsilon(u^\varepsilon)$.
\begin{longlist}[\textit{Step} 1.]
\item[\textit{Step} 1.] As $\int_\mathbb{R} xu^\varepsilon_0(x)\,dx=0$ and $\int_\mathbb{R} x^{8q^2}u^\varepsilon_0(x)\,dx<+\infty$, we have
\[
\Upsilon_\varepsilon \bigl(u^\varepsilon_0 \bigr)\geq
\inf_{u\in\mathcal{H}}  \Upsilon_\varepsilon(u)>\Upsilon_\varepsilon
(u_0 ) .
\]
We deduce $u^\varepsilon\neq u^\varepsilon_0$ since $t\mapsto\xi(t)=\Upsilon_\varepsilon (u_t^\varepsilon )$ is nonincreasing.

\item[\textit{Step} 2.] We proceed now a reductio ad absurdum by assuming
$u^\varepsilon=u^\varepsilon_-$. There exists $t_0>0$ such that $\int_\mathbb
{R} xu_{t_0}^\varepsilon(x)\,dx=0$. Consequently,
\[
\Upsilon_\varepsilon \bigl(u_{t_0}^\varepsilon \bigr)\geq
\inf_{u\in\mathcal
{H}} \Upsilon_\varepsilon(u)>\Upsilon_\varepsilon
(u_0 ),
\]
which contradicts the fact that $\xi$ is nonincreasing.

\item[\textit{Step} 3.] Assumption (M3) implies the weak convergence toward
$u^\varepsilon_+$.\quad\qed
\end{longlist}
\noqed\end{pf}
We use now Theorem~\ref{thmbaasym} in some particular cases.
%th7 #&#
%
\begin{thmm}
\label{thmbaasymapplication}
Let $du_0$ be a probability measure which verifies \textup{(FE)} and \textup{(FM)}. We
assume also
\[
\Upsilon (u_0 )<V(x_0)+\frac{1}{4}F(2x_0)
\quad \mbox{and}\quad  \int_\mathbb{R} xu_0(x)\,dx>0,
\]
where $x_0$ is defined in the \hyperref[intr]{Introduction}. Under either
conditions~\textup{(LIN)} or \textup{(SYN)}, $u_t^\varepsilon$ converges weakly toward
$u^\varepsilon_+$ for $\varepsilon$ small enough.
\end{thmm}
\begin{pf}
\textit{Step} 1. Theorem 3.2 in~\cite{HT1} and Theorem~\ref{trinity} imply
condition \textup{(M3)} under \textup{(LIN)} or \textup{(SYN)}.\vspace*{-6pt}
\begin{longlist}[\textit{Step} 2.]
\item[\textit{Step} 2.] Lemma~\ref{lemaenergylimit} provides the limit $\lim_{\varepsilon\longrightarrow0}\Upsilon_\varepsilon(u^\varepsilon_0)=V(x_0)+\frac
{1}{4}F(2x_0)$. Then, we deduce
%e9 #&#
%
\renewcommand{\theequation}{\arabic{section}.\arabic{equation}}
\setcounter{equation}{0}
\begin{equation}
\label{eqbaliminf} \lim_{\varepsilon\longrightarrow0}\inf_{u\in\mathcal{H}}
\Upsilon_\varepsilon (u)\leq V(x_0)+\frac{1}{4}F(2x_0)
 .
\end{equation}
\item[\textit{Step} 3.] We prove now that $  V(x_0)+\frac
{1}{4}F(2x_0)=\lim_{\varepsilon\longrightarrow0}\inf_{u\in\mathcal
{H}}\Upsilon_\varepsilon(u)$. Indeed, if $u$ is a probability measure such
that $\int_\mathbb{R} xu(x)\,dx=0$, it verifies the following inequality:
\begin{eqnarray*}
\Upsilon_\varepsilon(u)&\geq&\Upsilon_\varepsilon^-(u)+\frac{F''(0)}{4}
\iint_{\mathbb{R}^2}(x-y)^2u(x)u(y)\,dx\,dy
\\
&\geq&\Upsilon_\varepsilon^-(u)+\frac{F''(0)}{2}\int_\mathbb{R}
x^2u(x)\,dx .
\end{eqnarray*}
By using \eqref{eqfrminoration}, it yields
%e10 #&#
%
\begin{equation}
\label{eqbaminoration}   \Upsilon_\varepsilon(u)\geq-
\frac{\varepsilon}{4}-\frac{4\varepsilon
}{\exp(1)}+\int_\mathbb{R} \biggl
\{V(x)+\frac{F''(0)}{2}x^2-\frac
{\varepsilon x^2}{4} \biggr\}u(x)\,dx .
\end{equation}
We split now the study depending on whether we use conditions \textup{(LIN)} or~\textup{(SYN)}:
\begin{longlist}[(SYN)]
\item[(LIN)] If $F'$ is linear, $\frac{F''(0)}{2}x^2=\frac{1}{4}F(2x)$.
So the minimum of $x\mapsto V(x)+\frac{1}{4}F(2x)$ is $V(x_0)+\frac
{1}{4}F(2x_0)$. We can easily prove that
\[
\min_{x\in\mathbb{R}} \biggl(V(x)+\frac{F''(0)}{2}x^2-
\frac{\varepsilon
}{4}x^2 \biggr)=V(x_0)+
\frac{1}{4}F(2x_0)+o(1) .
\]
Consequently,
\[
\Upsilon_\varepsilon(u)\geq-\frac{\varepsilon}{4}-\frac{4\varepsilon}{\exp
(1)}+V(x_0)+
\frac{1}{4}F(2x_0)+o(1)
\]
for all $u\in\mathcal{H}$. Then, $ \lim_{\varepsilon
\longrightarrow0}\min_{u\in\mathcal{H}}\Upsilon_\varepsilon(u)\geq
V(x_0)+\frac{1}{4}F(2x_0)$. Inequality \eqref{eqbaliminf} provides
$ \lim_{\varepsilon\longrightarrow0}\inf_{u\in\mathcal
{H}}\Upsilon_\varepsilon(u)=V(x_0)+\frac{1}{4}F(2x_0)$.
\item[(SYN)] Since $V''(0)+F''(0)>0$, \eqref{eqbaminoration} implies
$\Upsilon_\varepsilon(u)\geq-\frac{\varepsilon}{4}-\frac{4\varepsilon}{\exp(1)}$
for all $u\in\mathcal{H}$ if $\varepsilon$ is less than $2
(V''(0)+F''(0) )$. We deduce\vadjust{\goodbreak} that\break $ \lim_{\varepsilon
\longrightarrow0}\inf_{u\in\mathcal{H}}\Upsilon_\varepsilon(u)\geq0$.
However, as $V''(0)+F''(0)>0$, Theorem 5.4 in~\cite{HT2} implies
$x_0=0$ so $V(x_0)+\frac{1}{4}F(2x_0)=0$. Inequality \eqref
{eqbaliminf} provides the following limit: $ \lim_{\varepsilon\longrightarrow0}\inf_{u\in\mathcal{H}}\Upsilon_\varepsilon
(u)=0=V(x_0)+\frac{1}{4}F(2x_0)$.
\end{longlist}
\item[\textit{Step} 4.] Consequently, $\Upsilon_\varepsilon(u_0)<\inf_{u\in\mathcal
{H}}\Upsilon_\varepsilon(u)$ for $\varepsilon$ small enough. Then, we apply
Theorem~\ref{thmbaasym}.\quad\qed
\end{longlist}
\noqed\end{pf}
%
%re3.1 #&#
%
\begin{rem}
\label{rembamiror}
We can replace $\int_\mathbb{R} xu_0(x)\,dx>0$ by $\int_\mathbb{R}
xu_0(x)\,dx<0$ in Theorems~\ref{thmbaasym} and \ref
{thmbaasymapplication}; then the same results hold with $u^\varepsilon_-$ instead of $u^\varepsilon_+$.
\end{rem}
%
%apA #&#
%
\begin{appendix}
%sB #&#
\section*{Appendix: Useful technical results}
In this annex, we present some results used previously in the proofs of
the main theorems.

Proposition~\ref{propaothermeasure} allows us to ensure that even if
the free-energy does not reach its global minimum on the stationary
measure $u^\varepsilon_0$, if the unique symmetric stationary measure is
an adherence value, then it is unique.

Proposition~\ref{propafm} is a general result on the self-stabilizing
processes. Indeed, it is well known that $du_t^\varepsilon$ is absolutely
continuous with respect to the Lebesgue measure for all $t>0$.
Proposition~\ref{propafm} extends this instantaneous regularization
to the finiteness of all the moments.

Lemma~\ref{lemaenergylimit} consists in asymptotic computation of the
free-energy in the small-noise limit for some useful measures. Lemma
\ref{lemalaplace} use a Laplace method for making a tedious
computation which is necessary for avoiding to assume that each family
of stationary measures verify condition (H).

We present now the essential proposition for proving Theorem~\ref{convergence}.
%prB.1 #&#
%
\begin{propp}
\label{propaothermeasure}
Let $du_0$ be a probability measure which verifies \textup{(FE)} and~\textup{(FM)}. We assume
the existence of two polynomial functions $\mathcal{P}$ and $\mathcal{Q}$,
a smooth function $\varphi$ with compact support such that
$\llvert \varphi(x)\rrvert \leq\mathcal{P}(x)$ and
$\llvert \varphi'(x)\rrvert^2\leq\mathcal{Q}(x)$, $\kappa>0$ and two
sequences $(r_k)_{k}$ and $(s_k)_{k}$ which go to $\infty$ such that
for all $r_k\leq t\leq s_k<r_{k+1}$,
\begin{eqnarray*}
\kappa=\int_\mathbb{R}\varphi(x)u_{r_k}^\varepsilon(x)\,dx
\leq\int_\mathbb {R}\varphi(x)u_t^\varepsilon(x)\,dx
\leq\int_\mathbb{R}\varphi (x)u_{s_k}^\varepsilon(x)\,dx=2
\kappa .
\end{eqnarray*}
Then, there exists a stationary measure $u^\varepsilon$ which verifies
$\int_\mathbb{R}\varphi(x)u^\varepsilon(x)\,dx\in[\kappa;2\kappa]$ and an increasing
sequence $(q_k)_{k}$ which goes to $\infty$ such that $u_{q_k}^\varepsilon$
converges weakly toward $u^\varepsilon$.
\end{propp}
\begin{pf}
\textit{Step} 1. We will prove that
$ \liminf_{k\longrightarrow+\infty} (s_k-r_k )>0$. We
introduce the function
\[
\Phi(t):=\int_\mathbb{R}\varphi(x)u_t^\varepsilon(x)\,dx
 .
\]
This function is well defined since $\llvert \varphi\rrvert $ is bounded
by a
polynomial function. The derivation\vadjust{\goodbreak} of $\Phi$, the use of equation
\eqref{eqintropde} and an integration by parts lead to
\begin{eqnarray*}
\Phi'(t)&=&-\int_\mathbb{R}\varphi'(x)
\biggl\{\frac{\varepsilon}{2}\frac
{\partial}{\partial
x}u_t^\varepsilon(x)+u_t^\varepsilon(x)
\bigl(V'(x)+F'\ast u_t^\varepsilon (x)
\bigr) \biggr\}\,dx
\\
&=&-\int_\mathbb{R}\varphi'(x)
\eta_t(x)\,dx .
\end{eqnarray*}
The Cauchy--Schwarz inequality implies
\[
\bigl\llvert \Phi'(t)\bigr\rrvert \leq\sqrt{\bigl\llvert
\xi'(t)\bigr\rrvert }\sqrt{\int_\mathbb
{R} \bigl(\varphi'(x) \bigr)^2u_t^\varepsilon(x)\,dx},
\]
where we recall that $\xi(t)=\Upsilon_\varepsilon (u_t^\varepsilon
)$. The
function $ (\varphi' )^2$ is bounded by a polynomial function, and
$\int_\mathbb{R} x^{2N}u_t^\varepsilon(x)$ is uniformly bounded with
respect to $t\in\mathbb{R}_+$
for all $N\in\mathbb{N}$. So, there exists $C>0$ such that
$\int_\mathbb{R} (\varphi'(x) )^2u_t^\varepsilon(x)\,dx\leq C^2$
for all
$t\in\mathbb{R}_+$. We deduce
%e11 #&#
%
\renewcommand{\theequation}{\Alph{section}.\arabic{equation}}
\setcounter{equation}{0}
\begin{equation}
\label{eqamajoration}
 \bigl\llvert \Phi'(t)\bigr\rrvert \leq C
\sqrt{\bigl\llvert \xi'(t)\bigr\rrvert } .
\end{equation}
By definition of the two sequences $(r_k)_{k}$ and $(s_k)_{k}$, we have
\[
\Phi(s_k)-\Phi(r_k)=\kappa .
\]
Combining this identity with \eqref{eqamajoration}, it yields
\[
C\int_{r_k}^{s_k}\sqrt{\bigl\llvert
\xi'(t)\bigr\rrvert }\,dt\geq\kappa .
\]
We apply the Cauchy--Schwarz inequality, and we obtain
\[
C\sqrt{s_k-r_k}\sqrt{\xi(r_k)-
\xi(s_k)}\geq\kappa
\]
since $\xi$ is nonincreasing; see Proposition \ref
{propintrodecrease}. Moreover, $\xi(t)$ converges as $t$ goes to
$\infty$; see Lemma~\ref{lemfrconv}. It implies the convergence of
$\xi(r_k)-\xi(s_k)$ toward $0$ when $k$ goes to $+\infty$.
Consequently, $s_k-r_k$ converges toward $+\infty$ so $
\liminf_{k\longrightarrow+\infty}s_k-r_k>0$.\vspace*{-6pt}
\begin{longlist}[\textit{Step} 2.]
\item[\textit{Step} 2.] By Lemma~\ref{lemfrconv}, $\Upsilon_\varepsilon
(u_t^\varepsilon )-\Upsilon_\varepsilon(u^\varepsilon)=\int_t^\infty\xi'(s)\,ds$
converges toward $0$. As $\xi'$ is nonpositive, we deduce that
$\sum_{k=N}^\infty\int_{r_k}^{s_k}\xi'(s)\,ds$ converges also toward $0$
when $N$ goes to $+\infty$. As $ \liminf_{k\longrightarrow
+\infty}s_k-r_k>0$, we deduce that there exists an increasing sequence
$q_k\in[r_k;s_k]$ which goes to $\infty$ and such that $\xi'
(q_k )$ converges toward $0$ when $k$ goes to $\infty$.
Furthermore, $\int_\mathbb{R}\varphi(x)u_{q_k}^\varepsilon(x)\,dx\in[\kappa
;2\kappa]$ for all $k\in\mathbb{N}$.

\item[\textit{Step} 3.] By proceeding similarly as in the proof of Theorem \ref
{thmfrsubconv}, we extract a subsequence of $(q_k)_{k}$ (we continue
to write it $q_k$ for simplifying the reading) such that
$u_{q_k}^\varepsilon$ converges weakly toward a stationary measure
$u^\varepsilon$. Moreover, $u^\varepsilon$ verifies $\int_\mathbb{R}\varphi
(x)u^\varepsilon(x)\,dx\in[\kappa;2\kappa]$.\quad\qed
\end{longlist}
\noqed\end{pf}
We provide now a result which allows us to obtain the statements of the
main theorem (Theorem~\ref{convergence}) with a weaker condition:\vadjust{\goodbreak}
%prB.2 #&#
%
\begin{propp}
\label{propafm}
Let $du_0$ be a probability measure which verifies \textup{(FE)} and~\textup{(ES)}. Then,
for all $t>0$, $du_t^\varepsilon$ satisfies \textup{(FM)}.
\end{propp}
\begin{pf}
\textit{Step} 1. If $du_0$ verifies (FM), then $du_t^\varepsilon$ satisfies
(FM) for all $t>0$; see Theorem 2.13 in~\cite{HIP}. We assume now that
$du_0$ does not satisfy (FM). Let us introduce $l_0:=\min \{l\geq0
 \llvert   \mathbb{E} [X_0^{2l} ] \}=+\infty$. We
know that $\mathbb{E} [X_t^{2l_0-2} ]<+\infty$ for all $t\geq
0$.\vspace*{-6pt}
\begin{longlist}[\textit{Step} 2.]
\item[\textit{Step} 2.] Let $t_0>0$. We proceed a reductio ad absurdum by
assuming that $\mathbb{E} [X_{t_0}^{2l_0} ]=+\infty$. This
implies directly $\mathbb{E} [X_{t}^{2l_0} ]=+\infty$ for all
$t\in[0,t_0]$. We recall that $2m$ (resp., $2n$) is the degree of the
confining (resp., interaction) potential $V$ (resp., $F$). Also, $q:=\max
\{m  ;  n \}$. For all $t\in[0,t_0]$, the application
$x\mapsto F'\ast u_t^\varepsilon(x)$ is a polynomial function with
parameters $m_1(t),\ldots,m_{2n-1}(t)$, where $m_j(t)$ is the $j$th
moment of the law $du_t^\varepsilon$. We recall inequality \eqref{majoration},
\[
\sup_{1\leq j\leq8q^2}\sup_{t\in[0,t_0]}m_j(t)\leq M_0
 .
\]
Consequently, the application $x\mapsto V'(x)+F'\ast u_t^\varepsilon(x)$
is a polynomial function with degree $2q-1$. Furthermore, the principal
term does not depend of the moments of the law $du_t^\varepsilon$, so we
can write
\[
V'(x)+F'\ast u_t^\varepsilon(x)=
\kappa_{2q-1}x^{2q-1}+\mathcal{P}_t(x),
\]
where $\kappa_{2q-1}\in\mathbb{R}_+^*$ is a constant, and $\mathcal
{P}_t$ is a polynomial function with degree at most $2q-2$. Moreover,
$\mathcal{P}_t$ is parametrized by the $2n$ first moments only. Let
$l\in\mathbb{N}$. We introduce the function $\mathcal
{Q}_t(x):=2lx^{2q-1}\mathcal{P}_t(x)-l(2l-1)\varepsilon x^{2l-2}$. As
$\mathcal{Q}_t$ is a polynomial function of degree less than $2l+2q-3$,
we have the following inequality:
%e12 #&#
%
\renewcommand{\theequation}{\Alph{section}.\arabic{equation}}
\setcounter{equation}{1}
\begin{equation}
\label{eqapolynomial} 2l\kappa_{2q-1}x^{2l+2q-2}+
\mathcal{Q}_t(x)\geq C_l \bigl(x^{2l+2q-2}-1 \bigr),
\end{equation}
where $C_l$ is a positive constant. The application of Ito formula provides
\[
dX_t^{2l}=2lX_t^{2l-1}\sqrt{
\varepsilon}\,dB_t- \bigl[2l\kappa_{2q-1}X_t^{2l+2q-2}+
\mathcal{Q}_t (X_t ) \bigr]\,dt .
\]
After integration, we obtain
\begin{eqnarray*}
X_{t_0}^{2l}&=&X_0^{2l}+2l\sqrt{
\varepsilon}\int_0^{t_0}X_t^{2l-1}\,dB_t-
\int_0^{t_0} \bigl[2l\kappa_{2q-1}X_t^{2l+2q-2}+
\mathcal{Q}_t (X_t ) \bigr]\,dt
\\
&\leq&X_0^{2l}+2l\sqrt{\varepsilon}\int_0^{t_0}X_t^{2l-1}\,dB_t-
\int_0^{t_0}C_l \bigl(X_t^{2l+2q-2}-1
\bigr)\,dt
\end{eqnarray*}
after using \eqref{eqapolynomial}. We choose $l:=l_0+1-q$, and then
we take the expectation. We obtain
\[
0\leq\mathbb{E} \bigl[X_{t_0}^{2l_0+2-2q} \bigr]\leq
C_1-C_2\int_0^{t_0}
\mathbb{E} \bigl[X_t^{2l_0} \bigr]\,dt,
\]
where $C_1$ and $C_2$ are positive constants. Since $\mathbb{E}
[X_t^{2l_0} ]=+\infty$ for all $t\in[0;t_0]$, this contradicts the
inequality $0\leq\mathbb{E} [X_{t_0}^{2l_0+2-2q} ]$.
Consequently, for all $t_0>0$: $\mathbb{E} [X_{t_0}^{2l_0}
]<+\infty$.\vadjust{\goodbreak}

\item[\textit{Step} 3.] Let $T>0$ and $l_1\in\mathbb{N}$ such that $l_1\geq l_0$
where the integer $l_0$ is defined as previously: $l_0:=\min \{l\geq
0  \llvert   \mathbb{E} [X_0^{2l} ]=+\infty \}$.
If $l_1=l_0$, the application of \textit{Step} $2$ leads to $\mathbb{E}
[X_T^{2l_1} ]<+\infty$. If $l_1>l_0$, we put $t_i:=\frac
{i}{l_1+1-l_0}T$ for all $1\leq i\leq l_1+1-l_0$. We apply \textit{Step} $2$ to
$t_1$, and we deduce $\mathbb{E} [X_{t_1}^{2l_0} ]<+\infty$.
By recurrence, we deduce $\mathbb{E} [X_{t_i}^{2l_0+2i}
]<+\infty$ for all $1\leq i\leq l_1+1-l_0$, in particular $\mathbb
{E} [X_{t_{l_0-l_1}}^{2l_0+2(l_1-l_0)} ]<+\infty$ that means
$\mathbb{E} [X_T^{2l_1} ]<+\infty$. This inequality holds for
all $l_1\geq l_0$, so the probability measure $du_T^\varepsilon$ satisfies
(FM).\quad\qed
\end{longlist}
\noqed\end{pf}
In order to obtain the thirdness of the stationary measure (or a weaker
result, see Theorem~\ref{trinity}), we need to compute the small-noise
limit of the free-energy for the stationary measures $u^\varepsilon_+$,
$u^\varepsilon_-$ and $u^\varepsilon_0$.
%leB.3 #&#
%
\begin{lemm}
\label{lemaenergylimit}
Let $\varepsilon_0$ such that there exist three families of stationary
measures $ (u^\varepsilon_+ )_{\varepsilon\in]0;\varepsilon_0]}$, $
(u^\varepsilon_- )_{\varepsilon\in]0;\varepsilon_0]}$ and $ (u^\varepsilon_0 )_{\varepsilon\in]0;\varepsilon_0]}$ which verify
\[
\lim_{\varepsilon\to0}u^\varepsilon_\pm=\delta_{\pm a}
\qquad\mbox{and}\qquad \lim_{\varepsilon\to0}u^\varepsilon_0=
\frac{1}{2}\delta_{x_0}+\frac{1}{2}\delta_{-x_0},
\]
where $x_0$ is defined in the \hyperref[intr]{Introduction}. Then, we have the following limits:
\[
\lim_{\varepsilon\to0}\Upsilon_\varepsilon \bigl(u^\varepsilon_\pm
\bigr)=V(a)\qquad \mbox{and}\qquad \lim_{\varepsilon\to0}\Upsilon_\varepsilon
\bigl(u^\varepsilon_0 \bigr)=V(x_0)+
\frac{1}{4}F (2x_0 ) .
\]
Plus, by considering the measure $v^\varepsilon_+(x):=Z^{-1}\exp
[-\frac{2}{\varepsilon} (V(x)+F(x-a) ) ]$, we have
\[
\lim_{\varepsilon\to0}\Upsilon_\varepsilon \bigl(v^\varepsilon_+ \bigr)=V(a)
 .
\]
\end{lemm}
\begin{pf}
\textit{Step} 1. We begin to prove the result for $u^\varepsilon_0$.\vspace*{-6pt}
\begin{longlist}[\textit{Step} 1.1.]
\item[\textit{Step} 1.1.] We can write $u^\varepsilon_0(x)=Z^{-1} \exp [-\frac
{2}{\varepsilon} (V(x)+F\ast u^\varepsilon_0(x) ) ]$ since it
is a stationary measure. Hence
\begin{eqnarray*}
\Upsilon_\varepsilon \bigl(u^\varepsilon_0 \bigr)&=&-
\frac{\varepsilon}{2}\log \biggl(\int_\mathbb{R} \exp \biggl[-
\frac{2}{\varepsilon} \bigl(V(x)+F\ast u^\varepsilon_0(x) \bigr)
\biggr]\,dx \biggr)
\\
&&{}-\frac{1}{2}\iint_{\mathbb{R}^2}F(x-y)u^\varepsilon_0(x)u^\varepsilon_0(y)\,dx\,dy
 .
\end{eqnarray*}
It has been proved in~\cite{HT3} [Theorem 1.2 if $V''(0)+F''(0)>0$,
Theorem 1.4 if $V''(0)+F''(0)=0$ and Theorem 1.3 if $V''(0)+F''(0)<0$
applied with $f_{2l}(x):=x^{2l}$] that the $2l$th moment of $u^\varepsilon_0$
tends toward $x_0^{2l}$ for all $l\in\mathbb{N}$. Since $F$ is a
polynomial function, we deduce the convergence of $\iint_{\mathbb
{R}^2}F(x-y)u^\varepsilon_0(x)u^\varepsilon_0(y)\,dx\,dy$ toward $\frac
{F(2x_0)}{2}$.

\item[\textit{Step} 1.2.] If $V''(0)+F''(0)\neq0$, we can apply Lemma A.4 in \cite
{HT3} to $f(x):=1$ and $U_\varepsilon(x):=V(x)+F\ast u^\varepsilon_0(x)$.
This provides
\[
\int_\mathbb{R} \exp \biggl[-\frac{2}{\varepsilon} \bigl(V(x)+F
\ast u^\varepsilon_0(x) \bigr) \biggr]\,dx=C_\varepsilon\exp
\biggl[-\frac{2}{\varepsilon} \biggl(V(x_0)+\frac{F(2x_0)}{2} \biggr)
\biggr],\vadjust{\goodbreak}
\]
where the constant $C_\varepsilon$ verifies $\varepsilon\log(C_\varepsilon
)\longrightarrow0$ in the small-noise limit. We deduce
\[
-\frac{\varepsilon}{2}\log \biggl(\int_\mathbb{R} \exp \biggl[-
\frac
{2}{\varepsilon} \bigl(V(x)+F\ast u^\varepsilon_0(x) \bigr)
\biggr]\,dx \biggr)\longrightarrow V(x_0)+\frac{F(2x_0)}{2}
\]
when $\varepsilon$ collapses. Consequently, it leads to the following limit:
\[
\Upsilon_\varepsilon \bigl(u^\varepsilon_0 \bigr)
\longrightarrow V(x_0)+\tfrac
{1}{4}F(2x_0) .
\]
\item[\textit{Step} 1.3.] We assume now $V''(0)+F''(0)=0$. Then $x_0=0$ according
to Proposition 3.7 and Remark 3.8 in~\cite{HT2}. Propositions~3.5 and
3.6 in~\cite{HT3} imply
\[
0<\liminf_{\varepsilon\to0}\varepsilon^{{1}/{(2m_0)}} \int_\mathbb{R}
\exp \biggl[-\frac{2}{\varepsilon} \bigl(V(x)+F\ast u^\varepsilon_0(x)
\bigr) \biggr]\,dx
\]
and
\[
\limsup_{\varepsilon\to0}\varepsilon^{{1}/{(2m_0)}} \int
_\mathbb{R} \exp \biggl[-\frac{2}{\varepsilon} \bigl(V(x)+F\ast
u^\varepsilon_0(x) \bigr) \biggr]\,dx<+\infty,
\]
where $m_0\in\mathbb{N}^*$ depends only on $V$ and $F$. We deduce
\[
-\frac{\varepsilon}{2}\log \biggl(\int_\mathbb{R} \exp \biggl[-
\frac
{2}{\varepsilon} \bigl(V(x)+F\ast u^\varepsilon_0(x) \bigr)
\biggr]\,dx \biggr)\longrightarrow0
\]
when $\varepsilon$ collapses. Consequently, we obtain the following limit:
\[
\Upsilon_\varepsilon \bigl(u^\varepsilon_0 \bigr)
\longrightarrow0=V(x_0)+\tfrac
{1}{4}F(2x_0) .
\]
\item[\textit{Step} 2.] We prove now the result for $u^\varepsilon_+$ (the proof is
similar for $u^\varepsilon_-$).

\item[\textit{Step} 2.1.] We can write $u^\varepsilon_+(x)=Z^{-1} \exp [-\frac
{2}{\varepsilon} (V(x)+F\ast u^\varepsilon_+(x) ) ]$ since it
is a stationary measure. Hence
\begin{eqnarray*}
\Upsilon_\varepsilon \bigl(u^\varepsilon_+ \bigr)&=&-\frac{\varepsilon}{2}\log
\biggl(\int_\mathbb{R} \exp \biggl[-\frac{2}{\varepsilon}
\bigl(V(x)+F\ast u^\varepsilon_+(x) \bigr) \biggr]\,dx \biggr)
\\
&&{}-\frac{1}{2}\iint_{\mathbb{R}^2}F(x-y)u^\varepsilon_+(x)u^\varepsilon_+(y)\,dx\,dy
 .
\end{eqnarray*}
It has been proved in~\cite{HT3} (Theorem 1.5 applied with
$f_l(x):=x^l$) that the $l$th moment of $u^\varepsilon_+$ tends toward
$a^l$ for all $l\in\mathbb{N}$. Since $F$ is a polynomial function, we
obtain the convergence of $\iint_{\mathbb{R}^2}F(x-y)u^\varepsilon_+(x)u^\varepsilon_+(y)\,dx\,dy$ toward $0$.

\item[\textit{Step} 2.2.] Since the second derivative of the application
$x\mapsto V(x)+F(x-a)$ in $a$ is positive, we can apply Lemma A.4 in
\cite{HT3} to $f(x):=1$ and $U_\varepsilon(x):=V(x)+F\ast u^\varepsilon_+(x)$
[after noting that $U_\varepsilon^{(i)}(x)$ tends toward
$V^{(i)}(x)+F^{(i)}(x-a)$ uniformly on each compact for all $i\in\mathbb
{N}$]. This provides
\[
\int_\mathbb{R} \exp \biggl[-\frac{2}{\varepsilon} \bigl(V(x)+F
\ast u^\varepsilon_+(x) \bigr) \biggr]\,dx=C_\varepsilon\exp \biggl[-
\frac{2}{\varepsilon}V(a) \biggr],
\]
where the constant $C_\varepsilon$ verifies $\varepsilon\log(C_\varepsilon
)\longrightarrow0$ in the small-noise limit. We deduce
\[
-\frac{\varepsilon}{2}\log \biggl(\int_\mathbb{R} \exp \biggl[-
\frac
{2}{\varepsilon} \bigl(V(x)+F\ast u^\varepsilon_+(x) \bigr) \biggr]\,dx \biggr)
\longrightarrow V(a)
\]
when $\varepsilon\longrightarrow0$. Consequently, the following limit holds:
\[
\Upsilon_\varepsilon \bigl(u^\varepsilon_0 \bigr)
\longrightarrow V(a) .
\]
\item[\textit{Step} 3.] We proceed similarly for $v^\varepsilon_+$.\quad\qed
\end{longlist}
\noqed\end{pf}
We provide here a useful asymptotic result linked to the Laplace method.
%leB.4 #&#

\begin{lemm}
\label{lemalaplace}
Let $U_{k}$ and $U\in\mathcal{C}^{\infty} (\mathbb{R},\mathbb
{R} )$ such that for all $i\in\mathbb{N}$, $U_k^{(i)}$ converges
toward $U{(i)}$ uniformly on each compact subset when $k$ goes to
$+\infty$. Let $ (\varepsilon_k )_k$ be a sequence which
converges toward $0$ as $k$ goes to $+\infty$. If $U$ has $r$ global
minimum locations $A_1<\cdots<A_r$ and if there exist $R>0$ and $k_c$
such that $U_k(x)>x^2$ for all $|x|>R$ and $k>k_c$, then, for $k$ big
enough, we have:
\begin{longlist}[(1)]
\item[(1)] $U_k$ has exactly one global minimum location $A_j^{(k)}$ on
each interval $I_j$, where~$I_j$ represents the Vorono\"{i} cells corresponding to the central points $A_j$, with $1\leq j\leq
r$.

\item[(2)] $A_j^{(k)}$ tends toward $A_j$ when $k$ goes to $+\infty$.
\end{longlist}

Furthermore, for all $N\in\mathbb{N}$, there exist $p_1,\ldots,p_r$
which verify $p_1+\cdots+p_r=1$ and $p_i\geq0$ for all $1\leq i\leq r$
such that we can extract a subsequence $\psi(k)$ which satisfies
\[
\lim_{k\to+\infty}\frac{\int_\mathbb{R} x^{l}\exp [-
{2}/{\varepsilon_{\psi(k)}}U_{\psi(k)} ]\,dx}{\int_\mathbb{R}\exp
[-{2}/{\varepsilon_{\psi(k)}}U_{\psi(k)} ]\,dx}=\sum_{j=1}^rp_jA_j^l
\]
for all $1\leq l\leq N$.
\end{lemm}
\begin{pf}
(1) The first point of the lemma is exactly the one of Lem\-ma~A.4
in~\cite{HT3}.\vspace*{-6pt}
\begin{longlist}[(2)]
\item[(2)] Since $U_k(x)\geq x^2$ for $|x|\geq R$ and $k>k_c$, we can
confine\vspace*{1pt} each $A_j^{(k)}$ in a compact subset. Then, the uniform
convergence on all the compact subset implies the convergence of
$A_j^{(k)}$ toward $A_j$ when $k$ goes to $+\infty$.\vspace*{1pt}

\item[(3)] Let $\rho>0$ arbitrarily small such that $[A_j-\rho,A_j+\rho
]\subset I_j$. For obvious reasons, we can extract a subsequence such that
\begin{eqnarray*}
\frac{\int_{A_i-\rho}^{A_i+\rho}\exp [-{2}/{\varepsilon_{\psi
(k)}}U_{\psi(k)}(x) ]\,dx}{\sum_{j=1}^r\int_{A_j-\rho}^{A_j+\rho}\exp
[-{2}/{\varepsilon_{\psi(k)}}U_{\psi(k)}(x)
]\,dx}\longrightarrow\lambda_i(\rho)
\end{eqnarray*}
with $\lambda_i(\rho)\geq0$ for all $1\leq i\leq r$ and $\sum_{j=1}^r\lambda_j(\rho)=1$.
\end{longlist}

\textit{We can note that the generation of
the sequence $\psi(k)$ depends on the choice of $\rho$. Consequently,
in the following, we can take $\rho$ arbitrarily small, then $\varepsilon_{\psi(k)}$ arbitrarily small.}

As the $r$ families $ (\lambda_j(\rho) )_{\rho>0}$ are
bounded, we can extract a subsequence $ (\rho_p )_{p}$ such
that $\lambda_j(\rho_p)$ tends toward $\lambda_j$ when $p$ goes to
$+\infty$. Furthermore, $\lambda_j\geq0$ for all $1\leq j\leq r$ and
$\sum_{j=1}^r\lambda_j=1$. For simplifying, we will write $\rho$
(resp.,
$k$) instead of $\rho_p$ [resp., $\psi(k)$].

We introduce the function $\zeta_l^{(k)}(x):=x^{l}\exp [-\frac
{2}{\varepsilon_k}U_{\varepsilon_k}(x) ]$ for all $l\in\mathbb{N}$. By
using classical analysis' inequality, we obtain
%e13 #&#
%
\renewcommand{\theequation}{\Alph{section}.\arabic{equation}}
\setcounter{equation}{2}
\begin{equation}
\label{eqabig} \Biggl\llvert \frac{\int_\mathbb{R}\zeta_l^{(k)}(x)\,dx}{\int_\mathbb{R}\zeta_0^{(k)}(x)\,dx}-\sum
_{j=1}^r\lambda_jA_j^l
\Biggr\rrvert \leq\mathcal {T}_1+\mathcal{T}_2+
\mathcal{T}_3+\mathcal{T}_4+\mathcal{T}_5
\end{equation}
with
\begin{eqnarray*}
\mathcal{T}_1(\rho)&:=&\Biggl\llvert \sum
_{j=1}^r \bigl(\lambda_j-
\lambda_j(\rho) \bigr)A_j^l\Biggr\rrvert
,\qquad \mathcal{T}_2 (\rho,R ):=\rho lR^{l-1}\mbox{,}
\\
\mathcal{T}_3 (\rho,k )&:=&\sum_{j=1}^r
\frac{\int_{I_j\cap
[A_j-\rho,A_j+\rho ]^c}\zeta_l^{(k)}(x)\,dx}{\int_\mathbb
{R}\zeta_0^{(k)}(x)\,dx},\qquad  \mathcal{T}_4 (R,k ):=2
\frac
{\int_R^{+\infty}\zeta_l^{(k)}(x)\,dx}{\int_\mathbb{R}\zeta_0^{(k)}(x)\,dx}
\end{eqnarray*}
and
\[
\mathcal{T}_5 (\rho,R,k ):=\sum_{j=1}^r
\llvert A_j\rrvert^l\biggl\llvert \frac{\int_{A_j-\rho}^{A_j+\rho}\zeta_0^{(k)}(x)\,dx}{\int_\mathbb
{R}\zeta_0^{(k)}(x)\,dx}-
\lambda_j(\rho)\biggr\rrvert \leq \Biggl(\sum
_{j=1}^r\llvert A_j
\rrvert^l \Biggr) (\mathcal{T}_3+\mathcal{T}_4
) .
\]
Let $\tau>0$ arbitrarily small. We take $R\geq2$ such that
\[
\max_{z\in[A_1-1;A_1+1]}U(z)+2<\frac{R^2}{2}.
\]
\begin{longlist}[3.1.]
\item[3.1.] The convergence of $\lambda_j(\rho)$ toward $\lambda_j$
implies the existence of $\rho_0>0$ such that for all $\rho<\rho_0$, we have
%e14 #&#
%
\begin{equation}
\label{eqasmall1} \mathcal{T}_1(\rho)\leq\frac{\tau}{5}
\end{equation}
for all $1\leq l\leq N$.

\item[3.2.] By taking $\rho<\min \{\rho_0 ; \min_{1\leq l\leq
N}\frac{\tau}{5lR^{l-1}} \}$, we deduce
%e15 #&#
%
\begin{equation}
\label{eqasmall2} \mathcal{T}_2 (\rho,R )\leq\frac{\tau}{5}
\end{equation}
for all $1\leq l\leq N$.

\item[3.3.] We will prove that the third term tends toward $0$. It is
sufficient to prove the following convergence:
\[
\frac{\int_{I_j\cap [A_j-\rho,A_j+\rho ]^c}\zeta_l^{(k)}(x)\,dx}{\int_{ [A_j-\rho,A_j+\rho ]}\zeta_0^{(k)}(x)\,dx}\longrightarrow0
\]
for all $1\leq j\leq r$. Since $I_j\subset[-R,R]$, we have
\begin{eqnarray*}
\frac{\int_{I_j\cap [A_j-\rho,A_j+\rho ]^c}\zeta_l^{(k)}(x)\,dx}{\int_{ [A_j-\rho,A_j+\rho ]}\zeta_0^{(k)}(x)\,dx}&\leq&R^l\frac{\int_{I_j\cap [A_j-\rho,A_j+\rho
]^c}\zeta_0^{(k)}(x)\,dx}{\int_{ [A_j-\rho,A_j+\rho ]}\zeta_0^{(k)}(x)\,dx} .
\end{eqnarray*}
Let us prove the convergence toward $0$ of the right-hand term:
\begin{eqnarray*}
&&\frac{\int_{I_j\cap [A_j-\rho,A_j+\rho ]^c}\zeta_0^{(k)}(x)\,dx}{\int_{ [A_j-\rho,A_j+\rho ]}\zeta_0^{(k)}(x)\,dx}\\
&&\qquad\leq R^{l+1}\frac{\sup \{\zeta_0^{(k)}(z)  ;
z\in I_j\cap [A_j-\rho,A_j+\rho ]^c \}}{\int_{A_j-
{\rho}/{2}}^{A_j+{\rho}/{2}}\zeta_0^{(k)}(x)\,dx}
\\
&&\qquad\leq\frac{R^{l+1}}{\rho}\frac{\sup \{\zeta_0^{(k)}(z)  ;  z\in
I_j\cap [A_j-\rho,A_j+\rho ]^c \}}{\inf \{\zeta_0^{(k)}(z)  ;  z\in [A_j-{\rho}/{2},A_j+{\rho
}/{2} ] \}}
\\
&&\qquad\leq\frac{R^{l+1}}{\rho}\exp \biggl\{-\frac{2}{\varepsilon_k} \Bigl[
\inf_{z\in I_j\cap [A_j-\rho,A_j+\rho ]^c}U_k(z)-\sup_{z\in
[A_j-{\rho}/{2},A_j+{\rho}/{2} ]}U_k(z)
\Bigr] \biggr\}.
\end{eqnarray*}
Let $\rho_1>0$ such that for all $\rho<\rho_1$, we have
\begin{eqnarray*}
\min_{1\leq j\leq r} \Bigl\{\inf_{z\in I_j\cap [A_j-\rho,A_j+\rho
]^c}U(z)-\sup_{z\in [A_j-{\rho}/{2},A_j+{\rho
}/{2} ]}U(z)
\Bigr\}\geq\delta>0 .
\end{eqnarray*}
We take\vspace*{1.5pt} $\rho<\min \{\rho_0,\rho_1,\min_{1\leq l\leq N}\frac{\tau
}{5lR^{l-1}} \}$. As $U_k$ converges uniformly toward $U$ on all
the compact subset, we deduce that for $k\geq k_0$, we have
%e16 #&#
%
\begin{equation}
\label{eqasmall3} \mathcal{T}_3 (\rho,k )\leq\frac{\tau}{5 (1+\max_{1\leq
l\leq N}\sum_{j=1}^r\llvert A_j\rrvert^l )}
\end{equation}
for all $1\leq l\leq N$.

\item[3.4.] By using the growth property on $U_k$ then the change of
variable $x:=\sqrt{\varepsilon_k}y$, it yields
\begin{eqnarray*}
\int_R^{+\infty}\zeta_l^{(k)}(x)&
\leq&\int_R^{+\infty}x^l\exp \biggl[-
\frac{2}{\varepsilon_k}x^2 \biggr]\,dx\leq C(l)e^{-{R^2}/{\varepsilon
_k}}
\varepsilon_k^{{(l+1)}/{2}},
\end{eqnarray*}
where $C(l)$ is a constant. We recall the assumption
$ \max_{z\in[A_1-1;A_1+1]}U(z)+2<\frac{R^2}{2}$. Since $U_k$ converges toward
$U$ uniformly on each compact subset, we have\vadjust{\goodbreak} $\max_{z\in
[A_1-1;A_1+1]}U_k(z)+1<\frac{R^2}{2}$ for $k\geq k_1$ (independently of
$\rho$). Consequently,
\begin{eqnarray*}
\mathcal{T}_4 (R,k )&\leq&\frac{2C(l)\varepsilon_k^{
{(l+1)}/{2}}\exp [-{2}/{\varepsilon_k} (\max_{z\in
[A_1-1;A_1+1]}U_k(z)+1 ) ]}{\int_{A_1-1}^{A_1+1}\exp
[-{2}/{\varepsilon_k}U_k(z) ]\,dx}
\\
&\leq& C(l)\varepsilon_k^{{(l+1)}/{2}}\exp \biggl[-\frac{2}{\varepsilon_k}
\biggr] .
\end{eqnarray*}
For $k\geq k_2$, we have the inequality
\[
\varepsilon_k^{{(l+1)}/{2}}\exp \biggl[-\frac{2}{\varepsilon_k} \biggr]
\leq \frac{\tau}{5\max_{1\leq l\leq N}C(l)\times (1+\max_{1\leq l\leq
N}\sum_{j=1}^r\llvert A_j\rrvert^l )} .
\]
By taking $k\geq\max \{k_0,k_1,k_2 \}$, we obtain
%e17 #&#
%
\begin{equation}
\label{eqasmall4} \mathcal{T}_4 (R,k )\leq\frac{\tau}{5 (1+\max_{1\leq l\leq
N}\sum_{j=1}^r\llvert A_j\rrvert^l )}
\end{equation}
for all $1\leq l\leq N$.

\item[3.5.] By taking $\rho<\min \{\rho_0,\rho_1,\frac{\tau
}{5lR^{l-1}} \}$ and $k\geq\max \{k_0,k_1,k_2 \}$,
inequalities~\eqref{eqabig}--\eqref{eqasmall3} and \eqref{eqasmall4} provide
\[
\Biggl\llvert \frac{\int_\mathbb{R}\zeta_l^{(k)}(x)\,dx}{\int_\mathbb{R}\zeta_0^{(k)}(x)\,dx}-\sum_{j=1}^r
\lambda_jA_j^l\Biggr\rrvert <\tau
\]
for all $1\leq l\leq N$. This achieves the proof.\quad\qed
\end{longlist}
\noqed\end{pf}
%
%reB.5 #&#
%
\begin{remm}
\label{remanotweak}
This lemma seems weaker than Lemma A.4 in~\cite{HT3}. However, here, we
do not assume that the second derivative of $U$ is positive in all the
global minimum locations.
\end{remm}
\end{appendix}
\section*{Acknowledgment}
It is a great pleasure to thank Samuel Herrmann for his remarks
concerning this work.
Most of the ideas for this work were found while I was at the Institut
\'{E}lie Cartan in Nancy. And so, I want to mention that I would not
have been able to write this paper without the hospitality I received
from the beginning.

Finalement, un tr\`{e}s grand merci \`{a} Manue et \`{a} Sandra pour tout.

% imsref loaded by akundreckaite, 2012-07-18 14:43:03
%

%suskaldyti doi

\printaddresses

\end{document}